\documentclass[11pt]{article}

\usepackage{fullpage}
\usepackage{amsmath,amsthm,amssymb}
\usepackage{thmtools}
\usepackage{mathtools}
\usepackage[T1]{fontenc}
\usepackage[utf8]{inputenc}
\usepackage{graphicx}
\usepackage{enumitem}
\usepackage{authblk}
\usepackage{hyperref}
\usepackage{xcolor}
\usepackage{tikz}
\usepackage{lineno}

\usetikzlibrary{positioning,calc}
\usepackage{subcaption}

\newtheorem{theorem}{Theorem}[section]
\newtheorem{lemma}[theorem]{Lemma}
\newtheorem{corollary}[theorem]{Corollary}
\newtheorem{proposition}[theorem]{Proposition}
\theoremstyle{definition}

\newtheorem{remark}[theorem]{Remark}

\begin{document}

\title{Domination in token graphs of graphs with domination number two}

\author[1]{Germ\'an Ben\'itez-Bobadilla\thanks{\texttt{german@ciencias.unam.mx}}}
\author[1]{Fernando Esteban Contreras-Mendoza\thanks{\texttt{esteban.contreras@ciencias.unam.mx}}}
\author[2]{Sergio Gerardo G\'omez-Galicia\thanks{\texttt{sgomez@math.cinvestav.mx}}}
\author[3]{Ana Trujillo-Negrete\thanks{\texttt{atrujillo@im.unam.mx}}}

\affil[1]{Facultad de Ciencias, Universidad Nacional Aut\'onoma de M\'exico, Mexico}
\affil[2]{Departamento de Matem\'aticas, CINVESTAV, Mexico}
\affil[3]{Instituto de Matem\'aticas, Universidad Nacional Aut\'onoma de M\'exico, Mexico}

\date{}

\maketitle

%\tableofcontents
%\linenumbers

\begin{abstract}
    The $k$-token graph $F_k(G)$ of a graph $G$ is the graph whose vertices are
    the $k$-subsets of $V(G)$, with two $k$-subsets adjacent if and only if
    their symmetric difference is an edge of $G$. In this paper we study the
    domination number of $F_k(G)$ for graphs $G$ whose domination number
    $\gamma(G)$ is at least two, with a primary focus on the case $\gamma(G) =
    2$.  For $k=2$, we give an upper bound for $\gamma(F_2(G))$ in terms of a
    dominating pair of $G$, and show that this bound is sharp for every double star
    $S_{a,b}$ with $2\le a\le b$; moreover, we determine the exact value for every
    double star, namely $\gamma(F_2(S_{a,b}))=2a+b-1$ for $1\le a\le b$. We also
    prove that $\gamma(F_2(G))\ge\lceil n/2\rceil$ for every graph $G$ on $n$
    vertices with $\gamma(G)\geq2$, with equality for the graphs with the maximum
    number of edges among those with $\gamma(G)\geq2$.  Then, we focus on complete
    bipartite graphs; we give upper bounds for $\gamma(F_2(K_{r,t}))$ that match the
    exact values independently established in recent work, together with a
    short proof that $\gamma(F_2(K_{2,n-2}))=n-2$.  For $k\ge3$ we show that $\gamma(F_k(K_{2,n-2}))$
    equals, up to lower-order terms, the Tur\'an number $T(n-2,k,k-1)$, and we
    determine $\gamma(F_3(K_{2,n-2}))$ exactly, obtaining $T(n-2,3,2) +n -4 =
    n^2/4+O(n)$; for $k=4$, determining the leading constant is related to the
    still-open Tur\'an conjecture for 3-uniform hypergraphs.  Finally, for
    $n\geq 12$ even, we prove that $\gamma\left(F_3(K_{n/2,n/2})\right) =\frac{7}{48}n^2+O(n)$. 
\end{abstract}

% % % % % % % % % % % % % % % % % % % % % % % % % % % % % % % % % % % % % % % %
%  %  %  %  %  %  %  %  %  %  %  %  %  %  %  %  %  %  %  %  %  %  %  %  %  %  % 
%                                                                             %
%                             S  E  C  T  I  O  N                             %
%                                                                             %
%  %  %  %  %  %  %  %  %  %  %  %  %  %  %  %  %  %  %  %  %  %  %  %  %  %  % 
% % % % % % % % % % % % % % % % % % % % % % % % % % % % % % % % % % % % % % % %

\section{Introduction}

\begin{figure}[t]
	\centering
	\includegraphics[width=0.6\textwidth]{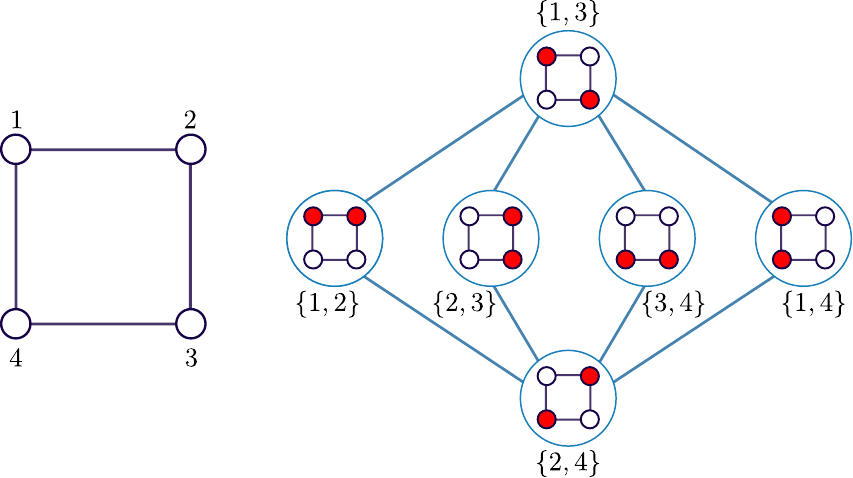}
	\caption{
		The cycle $C_4$ (left) and its $2$-token graph
		$F_2(C_4)\cong K_{2,4}$ (right).
		Each vertex of $F_2(C_4)$ represents a configuration
		of two indistinguishable tokens on $C_4$.
	}
	\label{fig:token-C4}
\end{figure}

Let $G$ be a graph on $n$ vertices and let $k\in\{1,\dots,n-1\}$. The
\emph{$k$-token graph} of $G$, denoted by $F_k(G)$, is the graph whose vertex set
is the family $\binom{V(G)}{k}$ of all $k$-subsets of $V(G)$, in which two
vertices $A$ and $B$ are adjacent if and only if their symmetric difference
$A\triangle B$ is an edge of $G$. Equivalently, the vertices of $F_k(G)$ are the
configurations of $k$ indistinguishable tokens placed on distinct vertices of
$G$, and two configurations are adjacent if and only if one can be obtained
from the other by moving a single token along an edge of $G$ to an unoccupied
vertex; for this reason, we also call the vertices of $F_k(G)$
\emph{configurations}.
It is easy to see that $F_1(G)\cong G$ and that the map $A\mapsto V(G)\setminus A$
is an isomorphism from $F_k(G)$ to $F_{n-k}(G)$~\cite{Token}.  Thus, it suffices
to consider $2\le k\le n/2$, and hence $n \ge 4$.  See
Figure~\ref{fig:token-C4} for the $2$-token graph of the cycle $C_4$.

Token graphs have been introduced independently several times and under
different names: as \emph{$k$-subgraph graphs} in Johns'
thesis~\cite{Johns}, as \emph{double vertex graphs} in the case
$k=2$~\cite{double_vertex}, as \emph{$n$-tuple vertex graphs}~\cite{ktuple}, and
as \emph{symmetric powers} of graphs~\cite{audenaert,rudolph}.
The name \emph{token graph}
and the systematic study of these graphs are due to Fabila-Monroy,
Flores-Pe\~naloza, Huemer, Hurtado, Urrutia, and Wood~\cite{Token}. Several
graph properties have been studied on token graphs so far, including connectivity~\cite{Token,connectivitytoken}, diameter, clique number,
chromatic number and Hamiltonicity~\cite{Token}, regularity and
planarity~\cite{regularity}, automorphisms~\cite{cuts_automorphisms},
independence and matching numbers~\cite{independence}, and Laplacian
spectra~\cite{Laplacian}.

A \emph{dominating set} of a graph $G$ is a set $D\subseteq V(G)$ such that
every vertex of $G$ either belongs to $D$ or has a neighbor in $D$, and the
\emph{domination number} $\gamma(G)$ of $G$ is the minimum size of a dominating
set of $G$. 
For background on domination in graphs we refer to~\cite{HHS98,HHH20}.
Since a vertex dominates at most $\Delta(G)+1$ vertices, where $\Delta(G)$
denotes the maximum degree of $G$, every graph $G$ on $n$ vertices satisfies
$\gamma(G)\geq n/(\Delta(G)+1)$~\cite{HHS98}; we use this bound in
Section~\ref{sec:bip3}. Note also that $\gamma(G)=1$ if and only if $G$ has a
\emph{universal} vertex, that is, a vertex adjacent to all the other vertices.

The domination number of token graphs has been studied only recently. Prabhu and
Manjula~\cite{PrabhuManjula2019} gave general bounds for $\gamma(F_2(G))$ in
terms of parameters of $G$. Fabila-Monroy and
G\'omez-Galicia~\cite{domination_token} gave asymptotic bounds for
$\gamma(F_k(K_{1,n}))$ and $\gamma(F_k(K_n))$, with the exact values when $k=2$,
namely $\gamma(F_2(K_{1,n}))=n-1$ and $\gamma(F_2(K_n))=\lfloor n/2\rfloor$.
The latter graphs are the Johnson
graphs, whose domination and total domination were also studied by Cornet,
Dravec, and Torres~\cite{CornetDravecTorres2026}, and the $k$-token graph of the star
$K_{1,n}$ for $k\geq3$ is studied by Acosta Troncoso,
Hern\'andez-Cruz, and the fourth author of this
paper~\cite{AcostaHernandezTrujillo}.
For paths, Acosta Troncoso, Lea\~nos, and
Rivera~\cite{troncoso} proved that $\gamma(F_2(P_n))=n^2/10+\Theta(n)$, and
Acosta Troncoso and Rivera~\cite{AcostaRivera2026} determined the exact value
for $n\ge13$. When this paper was in its final stage of preparation, two
preprints appeared, by Acosta Troncoso, Caldera Caldera, and
Rivera~\cite{AcostaCalderaRivera2026} and by Chinman, Cortijo, Pedroni, Seong,
and Stenger~\cite{Chinman2026}, which independently determine
$\gamma(F_2(K_{m,n}))$ for all $2\le m\le n$; the latter also gives bounds for
$\gamma(F_k(K_{m,n}))$ when $k\ge3$. Our results on $F_2(K_{r,t})$ in
Section~\ref{sec:bip2} were obtained independently of these works; we have kept
them because the constructions and the proof for $K_{2,n-2}$ are short and
self-contained.

In this paper we study the domination number of $F_k(G)$ when $\gamma(G)\ge2$,
that is, when no vertex of $G$ is adjacent to all the other vertices; note that
the star and the complete graph mentioned above are, respectively, the sparsest
and the densest graphs with $\gamma(G)=1$.

In Section~\ref{sec:gamma2} we consider the case $k=2$. We first give an upper
bound for $\gamma(F_2(G))$ in terms of a dominating pair of $G$, provided that
each of its two private neighborhoods has at least two vertices.

\begin{restatable}{theorem}{thmpair}
    \label{thm:pair}
    Let $G$ be a graph on $n$ vertices with $\gamma(G) = 2$, and let $\{x,y\}$
    be a dominating pair of $G$, that is, $N_G[x]\cup N_G[y]=V(G)$. Let $X
    \coloneqq N_G(x)\setminus N_G[y]$ and $Y  \coloneqq  N_G(y)\setminus
    N_G[x]$. If $|X|,|Y| \ge 2$, then
    \[
        \gamma\bigl(F_2(G)\bigr)\le n-3+\min\{|X|,|Y|\}.
    \]
\end{restatable}

For $a,b\ge1$, the \emph{double star} $S_{a,b}$ is the tree obtained from the
stars $K_{1,a}$ and $K_{1,b}$ by joining their centers with an edge; so
$S_{a,b}$ has $a+b+2$ vertices.  We show that the bound in Theorem
\ref{thm:pair} is attained by every double star $S_{a,b}$ with $2 \le a \le b$.
Moreover, we determine exactly the domination number of the $2$-token graph of
every double star.

\begin{restatable}{theorem}{thmdoublestars}\label{thm:doublestars}
    Let $1\leq a\leq b$ be integers. Then
    \[
        \gamma(F_2(S_{a,b}))=2a+b-1.
    \]
\end{restatable}

\noindent We then prove a general lower bound for $\gamma(F_2(G))$, which is
attained by the graphs with domination number two and the maximum possible
number of edges. 

\begin{restatable}{theorem}{thmdense}\label{thm:dense}
    Let $G$ be a graph on $n$ vertices with $\gamma(G)\geq2$. Then
    \[
        \gamma(F_2(G))\geq\left\lceil\frac{n}{2}\right\rceil,
    \]
    and equality holds if $G$ has the maximum number of edges among the graphs
    with domination number at least two, that is, if $\overline{G}\cong\frac
    n2K_2$ when $n$ is even, or $\overline{G}\cong P_3\cup\frac{n-3}2K_2$ when
    $n$ is odd.
\end{restatable}

\noindent In Section~\ref{sec:bipartite} we study complete bipartite graphs. For
$k=2$ we give dominating sets of $F_2(K_{r,t})$ whose sizes match the exact
values obtained very recently in~\cite{AcostaCalderaRivera2026,Chinman2026}, and we determine
$\gamma(F_2(K_{2,n-2}))$. For $k\ge3$ we show that computing
$\gamma(F_k(K_{2,n-2}))$ is essentially a Tur\'an problem.

Let $m\geq k\geq3$. A family $\mathcal{F}\subseteq\binom{[m]}{k-1}$ is a
\emph{Tur\'an $(m,k,k-1)$-family} if every $k$-subset of $[m]$ contains a member
of $\mathcal{F}$; the minimum size of such a family is the \emph{Tur\'an number}
$T(m,k,k-1)$, and the limit
\[
    t_k \coloneqq \lim_{m\to\infty}\frac{T(m,k,k-1)}{\binom{m}{k-1}}
\]
exists; see Section~\ref{sec:bip3} for this and the other facts on Tur\'an numbers that we
use. By Mantel's theorem~\cite{Mantel1907}, $t_3=\frac12$, but already $t_4$ is unknown:
Tur\'an~\cite{Turan1941} conjectured that $t_4=\frac49$, and the best known bounds are
$0.4422\leq t_4\leq\frac49$~\cite{JeongEtAl2026}.

In the following statement, $C(m,k-1,k-2)$ denotes the minimum size of a family
of $(k-1)$-subsets of $[m]$ such that every $(k-2)$-subset of $[m]$ is contained
in some member of the family, and $T(m,k-1,k-2)$ is defined as above with $k-1$
in place of $k$; for instance, $T(m,2,1)=m-1$.

\begin{restatable}{theorem}{thmturan}\label{thm:turan}
Let $k\geq3$ and $n\geq2k$. Then
\[
T(n-2,k,k-1)\leq\gamma(F_k(K_{2,n-2}))
 \leq T(n-2,k,k-1)+T(n-2,k-1,k-2)+C(n-2,k-1,k-2).
\]
Moreover, if $\mathcal{F}_{k-1}$ is a Tur\'an $(n-2,k,k-1)$-family and
$\mathcal{F}_{k-2}$ is a Tur\'an $(n-2,k-1,k-2)$-family such that every
$(k-2)$-subset of $[n-2]$ that is not in $\mathcal{F}_{k-2}$ is contained in
some member of $\mathcal{F}_{k-1}$, then
$\gamma(F_k(K_{2,n-2}))\leq|\mathcal{F}_{k-1}|+|\mathcal{F}_{k-2}|$; in
particular, if such families of sizes $T(n-2,k,k-1)$ and $T(n-2,k-1,k-2)$ exist,
then $\gamma(F_k(K_{2,n-2}))\leq T(n-2,k,k-1)+T(n-2,k-1,k-2)$.  Consequently,
for fixed $k$,
\[
    \gamma\bigl(F_k(K_{2,n-2})\bigr) =
        T(n-2,k,k-1)+O(n^{k-2}) =
        \bigl(t_k+o(1)\bigr)\binom{n}{k-1},
\]
% where the term $O(n^{k-2})$ is at most
where the additional term is bounded above by
$\bigl(t_{k-1}+\frac{1}{k-1}+o(1)\bigr)\binom{n-2}{k-2}$.
\end{restatable}

\noindent For $k=3$ we determine $\gamma(F_3(K_{2,n-2}))$ exactly
(Theorem~\ref{thm:K2n3}): it equals $T(n-2,3,2)+n-4=\frac{n^2}{4}+O(n)$.  For
$k=4$, Tur\'an's conjecture implies that the leading constant is $\frac49$.
Finally, we consider the $3$-token graph of the balanced complete bipartite graph
$K_{n/2,n/2}$ for $n$ even, for which we prove the following.

\begin{restatable}{theorem}{thmKmm}\label{thm:Kmm}

Let $n\geq 12$ be even. Then
\[
\frac{7}{48}n^2-4n
\leq
\gamma\left(F_3(K_{n/2,n/2})\right)
\leq
\frac{7}{48}n^2+O(n).
\]
In particular,
\[
\gamma\left(F_3(K_{n/2,n/2})\right)
=\frac{7}{48}n^2+O(n).
\]
\end{restatable}

\paragraph{Notation.} All graphs in this work are simple graphs, i.e.\ undirected graphs
without loops or parallel edges. We denote by $V$ and $E$ the vertex set and the edge set of a
graph $G$, and we usually use $n$ and $m$ for $|V|$ and $|E|$, respectively. The
number of triangles of $G$ is denoted by $\operatorname{tr}(G)$. We denote by
$\overline{G}$ the complement of $G$ and write $uv$ for the edge $\{u,v\}$; for
two vertices $A$ and $B$ of a graph, we write $A\sim B$ when they are adjacent.
We use $N_G(u)$ and $N_G[u]$ to denote the neighborhood and the closed
neighborhood of the vertex $u$ in $G$, respectively. For a finite set $S$, we
denote by $\binom{S}{k}$ the family of all $k$-subsets of $S$, and we write
$[n]\coloneqq\{1,\dots,n\}$.

% % % % % % % % % % % % % % % % % % % % % % % % % % % % % % % % % % % % % % % %
%                                                                             %
%                             S  E  C  T  I  O  N                             %
%                                                                             %
% % % % % % % % % % % % % % % % % % % % % % % % % % % % % % % % % % % % % % % %

% \section{Graphs with domination number two}\label{sec:gamma2}
\section[Bounds on gamma(F2(G)) for graphs G with gamma(G) >= 2]
{Bounds on $\gamma(F_2(G))$ for graphs $G$ with \texorpdfstring{$\gamma(G)\ge
2$}{$\gamma(G)\geq 2$}}
\label{sec:gamma2}

Throughout this section we study the domination number of $F_2(G)$. We order the results by the number of edges of $G$: we start with an upper
bound in terms of a dominating pair of $G$, which is attained by double stars, and
end with a lower bound valid whenever $\gamma(G)\geq2$, which is attained by the
graphs with the maximum number of edges.

%  %  %  %  %  %  %  %  %  %  %  %  %  %  %  %  %  %  %  %  %  %  %  %  %  %  % 
%                         S  U  B  S  E  C  T  I  O  N                        %
%  %  %  %  %  %  %  %  %  %  %  %  %  %  %  %  %  %  %  %  %  %  %  %  %  %  % 

\subsection{A general upper bound}

Our first result establishes an upper bound for $\gamma(F_2(G))$ in terms of the
sizes of the two private neighborhoods of a dominating pair of $G$.

\thmpair*

\begin{proof}
    Let $Z \coloneqq N_G(x)\cap N_G(y)$, and set $a \coloneqq |X|$, $b \coloneqq
    |Y|$ and $c \coloneqq |Z|$. Since $\{x,y\}$ is a dominating pair, $V(G)$ is
    the disjoint union of $\{x,y\}$, $X$, $Y$ and $Z$, so $n=a+b+c+2$. By
    symmetry, we may assume that $a\leq b$, so we just need to show that
    $\gamma(F_2(G))\leq n-3+a$. Let $X=\{x_1,\dots,x_a\}$, $Y=\{y_1,\dots,y_b\}$
    and $Z=\{z_1,\dots,z_c\}$, and define
    \begin{align*}
        D\coloneqq &\underbrace{\bigl\{\{x,x_i\}:i\in [a-2]\bigr\}}_{D_X}
        \cup\underbrace{\bigl\{\{y,y_j\}:j\in [b-2]\bigr\}}_{D_Y}
        \cup\underbrace{\bigl\{\{y,z_p\}:p\in [c]\bigr\}}_{D_Z}
        \\[-2mm]
        &\cup\underbrace{\bigl\{\{y,x_q\}:q\in [a]\bigr\}}_{C_{XY}} \cup
        \underbrace{\bigl\{\{x,y\},\{x_{a-1},x_a\},\{y_{b-1},y_b\}\bigr\}}_{\Delta}.
    \end{align*}

    We claim that $D$ dominates $F_2(G)$. First, notice that
    \[
        |D|=(a-2)+(b-2)+c+a+3 = 2a+b+c-1 = (a+b+c+2)-3+a = n-3+a.
    \]

    Let $\{u,v\}\in V(F_2(G))$. We distinguish the following cases.

    \begin{enumerate}[label=(\arabic*)]
        \item \emph{Suppose $u=x$.} If $v=x_i$ with $1\le i\le a-2$ then
        $\{x,v\}\in D_X\subseteq D$.  If $v\in\{x_{a-1},x_a\}$ then $\{x,v\}$ is
        adjacent to $\{x_{a-1},x_a\}\in\Delta$.  If $v\in Y\cup Z$, then
        $\{x,v\}$ is adjacent to $\{x,y\}\in\Delta$.  If $v=y$, then
        $\{x,y\}\in\Delta$.

        \item \emph{Suppose $u=y$.} If $v=y_j$ with $1\le j\le b-2$ then
        $\{y,v\}\in D_Y\subseteq D$.  If $v\in\{y_{b-1},y_b\}$ then $\{y,v\}$ is
        adjacent to $\{y_{b-1},y_b\}\in\Delta$.  If $v\in X\cup Z$ then
        $\{y,v\}\in C_{XY}\cup D_Z\subseteq D$.

        \item \emph{Suppose $u,v\in X\cup Y\cup Z$.} If
        $\{u,v\}=\{y_{b-1},y_b\}$ we are done since it belongs to $\Delta$.
        Otherwise, choose $u\notin\{y_{b-1},y_b\}$ (possible since
        $\{u,v\}\neq\{y_{b-1},y_b\}$).  If $v$ is adjacent to $y$ (i.e., $v\in
        Y\cup Z$), then $\{u,v\}$ is adjacent to $\{y,u\}\in D_Y\cup D_Z\cup
        C_{XY}$, because $\{y,u\}\in D_Y$ if $u\in Y$ with index $\le b-2$, it
        equals $\{y,u\}\in D_Z$ if $u\in Z$, and it lies in $C_{XY}$ if $u\in
        X$.  If $v$ is not adjacent to $y$, then $v\in X$. If $u$ is adjacent to
        $y$, by symmetry $\{u,v\}$ is adjacent to $\{y,v\}\in C_{XY}$.  If $u$
        is also in $X$, then either $\{u,v\}=\{x_{a-1},x_a\}\in\Delta$, or,
        assuming $u\notin\{x_{a-1},x_a\}$, we have that $\{u,v\}$ is adjacent to
        $\{x,u\}\in D_X$.
    \end{enumerate}
    This exhausts all possibilities, so $D$ is a dominating set and
    \[
        \gamma(F_2(G))\le|D|=n-3+a.
    \]
\end{proof}

%  %  %  %  %  %  %  %  %  %  %  %  %  %  %  %  %  %  %  %  %  %  %  %  %  %  % 
%                         S  U  B  S  E  C  T  I  O  N                        %
%  %  %  %  %  %  %  %  %  %  %  %  %  %  %  %  %  %  %  %  %  %  %  %  %  %  % 

\subsection{Double stars}\label{sec:doublestars}

\begin{figure}[t]
\centering
\includegraphics[width=0.9\textwidth]{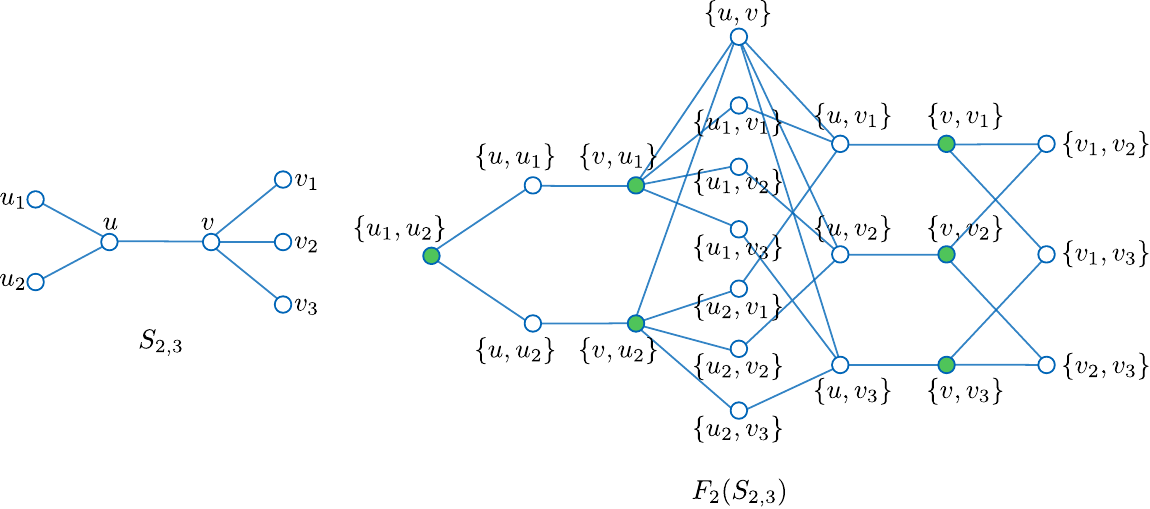}
\caption{The double star $S_{2,3}$ and its $2$-token graph $F_2(S_{2,3})$, drawn
according to the types of configurations used in the proof of
Theorem~\ref{thm:doublestars}. The filled vertices form the dominating set $D$
of that proof, which has $2a+b-1=6$ elements.}
\label{fig:doublestar}
\end{figure}

Among connected graphs with domination number two, those with the minimum number
of edges are trees, and the trees in which two adjacent vertices form a
dominating set are precisely the double stars. We now determine the domination number of their $2$-token graphs; see
Figure~\ref{fig:doublestar} for an example.

\thmdoublestars*

\begin{proof}
    Let $u$ and $v$ be the centers of $S_{a,b}$, with $uv\in E(S_{a,b})$, and
    let
    \[
        X=\{u_1,\ldots,u_a\},
        \qquad
        Y=\{v_1,\ldots,v_b\}
    \]
    be their respective sets of leaves. For convenience, define the following
    vertices of $F_2(S_{a,b})$:
    \begin{gather*}
        A_i=\{u,u_i\}, \qquad 
        B_j=\{v,v_j\}, \qquad 
        C_i=\{v,u_i\}, \qquad 
        E_j=\{u,v_j\},\\
        P_{ii'}=\{u_i,u_{i'}\},\qquad 
        Q_{jj'}=\{v_j,v_{j'}\}, \qquad 
        R_{ij}=\{u_i,v_j\}.
    \end{gather*}
    Figure~\ref{fig:doublestar} shows $F_2(S_{2,3})$ with its vertices arranged
    by these types. Note that
    \begin{gather*}
        N[P_{ii'}]=\{P_{ii'},A_i,A_{i'}\},\qquad 
        N[Q_{jj'}]=\{Q_{jj'},B_j,B_{j'}\},\qquad 
        N[R_{ij}]=\{R_{ij},C_i,E_j\},\\
        N[A_i]=\{A_i,C_i\}\cup\{P_{ii'}:i'\neq i\},\qquad 
        N[B_j]=\{B_j,E_j\}\cup\{Q_{jj'}:j'\neq j\},\\
        N[C_i]=\{C_i,A_i,\{u,v\}\}\cup\{R_{ij}:j\in[b]\},\qquad 
        N[E_j]=\{E_j,B_j,\{u,v\}\}\cup\{R_{ij}:i\in[a]\}.
    \end{gather*}

    We first establish the upper bound. Consider the following set
    \[
    D =
    \begin{cases}
        \{C_1\} \cup \{B_j \colon j\in[b]\}, & \text{if } a=1,\\
        \{C_i \colon i\in[a]\} \cup \{B_j \colon j\in[b]\} \cup
        \{A_i \colon i\in[a-2]\} \cup \{P_{a-1,a}\}, & \text{if }a\ge2.
    \end{cases}
    \]
    Hence, $|D| = 2a+b-1$.
    % \[
    %     D=\{C_i:i\in[a]\}
    %     \cup\{B_j:j\in[b]\}
    %     \cup\{A_i:i\in[a-2]\}
    %     \cup\{P_{a-1,a}\},
    % \]
    % where the last two sets are empty if $a=1$. Hence,
    % \[
    %     |D|=\begin{cases}
    %         a+b+(a-2)+1=2a+b-1 & \text{if $a\ge 2$,}\\
    %         1+b=2a+b-1 & \text{if $a=1$.}
    %     \end{cases}
    % \]

    We claim that $D$ dominates $F_2(S_{a,b})$.  Indeed, the configuration
    $\{u,v\}$ is adjacent to every $C_i$, and each $A_i$ is adjacent to $C_i$.
    Moreover, every $E_j$ and every $Q_{jj'}$ is adjacent to some $B_j\in D$,
    while every $R_{ij}$ is adjacent to $C_i\in D$. Finally, every $P_{ii'}$ is
    adjacent to some $A_i\in D$, unless $\{i,i'\}=\{a-1,a\}$, in which case
    $P_{ii'}\in D$. Therefore, $D$ is a dominating set, and
    \[
        \gamma(F_2(S_{a,b}))\leq 2a+b-1.
    \]

    We now prove the lower bound, and we consider first the case $a=1$. Let $D'$
    be an arbitrary dominating set of $F_2(S_{1,b})$. Suppose first that
    $C_1\notin D'$. Since $N[R_{1j}]=\{R_{1j},C_1,E_j\}$, at least one of
    $R_{1j}$ and $E_j$ belongs to $D'$ for each $j\in[b]$; these configurations
    are distinct for different values of $j$, so $D'$ contains at least $b$ of
    them. Moreover, $N[A_1]=\{A_1,C_1\}$, so $A_1\in D'$, and hence $|D'|\geq
    b+1$. Now suppose that $C_1\in D'$ and $\{u,v\}\notin D'$. Then, for every
    $j\in[b]$, the configuration $E_j$ must be dominated by an element of the
    set $\{E_j,B_j,R_{1j}\}$; these sets are pairwise disjoint as $j$ varies,
    and none of them contains $C_1$. Therefore, $|D'|\geq 1+b$. Finally, suppose
    that both $C_1$ and $\{u,v\}$ belong to $D'$. Let
    \[
        I=\{j\in[b]:B_j\in D'\},\qquad s=b-|I|.
    \]
    Since $N[Q_{jj'}]=\{Q_{jj'},B_j,B_{j'}\}$, every configuration $Q_{jj'}$
    with $j,j'\notin I$ can only be dominated by itself, so all such
    configurations belong to $D'$. Consequently,
    \[
        |D'|\geq 2+|I|+\binom{s}{2}=2+b-s+\binom s2=b+2+\frac{s(s-3)}{2}\geq b+1,
    \]
    where the last inequality holds for every integer $s\geq0$. Thus, in every
    case, $\gamma(F_2(S_{1,b}))\geq b+1=2a+b-1$.

    Suppose now that $2\leq a\leq b$. Let $D'$ be an arbitrary dominating set of
    $F_2(S_{a,b})$, and define the following sets
    \[
        I_A=\{i:A_i\in D'\}, \quad
        I_B=\{j:B_j\in D'\}, \quad
        I_C=\{i:C_i\in D'\}, \quad
        I_E=\{j:E_j\in D'\}.
    \]
    Let $p,q,r$ denote the number of configurations of types $P$, $Q$ and $R$,
    respectively, that belong to $D'$. Let
    \[
        z=
        \begin{cases}
        1,&\text{if }\{u,v\}\in D',\\
        0,&\text{otherwise}.
        \end{cases}
    \]

    Since $N[P_{ii'}]=\{P_{ii'},A_i,A_{i'}\}$, every pair of indices outside
    $I_A$ forces its corresponding configuration $P_{ii'}$ to belong to $D'$.
    Setting $s=a-|I_A|$, we obtain $p\geq\binom{s}{2}$, and consequently
    \[
        |I_A|+p
        \geq a-s+\binom{s}{2}
        =a+\frac{s(s-3)}2
        \geq a-1.
    \]
    Analogously,
    \[
        |I_B|+q\geq b-1.
    \]

    Now observe that $N[R_{ij}]=\{R_{ij},C_i,E_j\}$.  Set
    \[
        a'=a-|I_C|,
        \qquad
        b'=b-|I_E|.
    \]
    Every pair $(i,j)$ with $i\notin I_C$ and $j\notin I_E$ forces $R_{ij}\in
    D'$.  Thus $r\geq a'b'$, and hence
    \[
        |I_C|+|I_E|+r
        \geq a+b-a'-b'+a'b'\geq \min\{a,b\}=a,
    \]
    where the last inequality follows from
    \[
        a+b-a'-b'+a'b'-a = b-a'+(a'-1)b' \ge 0.
    \]

    Indeed, if $a'=0$, this expression is $b-b' \ge 0$; if $a'=1$, it equals
    $b-1>0$; and if $a' \ge 2$, both terms are nonnegative, since $a' \le a \le
    b$.  Moreover, equality occurs precisely when $(a',b') = (0,b)$, or, if
    $a=b$, when $(a',b') = (a,0)$.

    Combining these inequalities gives
    \[
        |D'|
        =(|I_A|+p)+(|I_B|+q)+(|I_C|+|I_E|+r)+z \geq (a-1)+(b-1)+a =2a+b-2.
    \]

    Suppose, aiming for a contradiction, that $|D'|=2a+b-2$. Then equality must
    hold in all three estimates, and $z=0$. In particular,
    \[
        |I_A|+p=a-1,\qquad
        |I_B|+q=b-1, \qquad \text{and}\qquad
        |I_C|+|I_E|+r=a.
    \]

    Since $2\leq a\leq b$, equality in the last inequality implies
    \[
        I_C=[a],\qquad I_E=\varnothing,\qquad r=0,
    \]
    or, when $a=b$, the symmetric possibility
    \[
        I_C=\varnothing,\qquad I_E=[b],\qquad r=0.
    \]

    In the first case, since $z=0$, $I_E=\varnothing$, and $r=0$, every
    configuration $E_j$ must be dominated by $B_j$. Thus $I_B=[b]$, which
    implies
    \[
        |I_B|+q\geq b,
    \]
    a contradiction.

    In the second case, we have $I_C=\varnothing$, $r=0$, and $z=0$. Since
    \[
        N[C_i]=\{C_i,A_i,\{u,v\}\}
        \cup\{R_{ij}:j\in[b]\},
    \]
    every configuration $C_i$ must be dominated by $A_i$. Consequently,
    $I_A=[a]$, which contradicts $|I_A|+p=a-1$.  Therefore,
    \[
        |D'|\geq 2a+b-1.
    \]
    Combining this with the upper bound yields the result.
\end{proof}

A double star $S_{a,b}$ with $2\leq a\leq b$ has $a+b+2$ vertices and its
centers form a dominating pair whose sets $X$ and $Y$, as defined in
Theorem~\ref{thm:pair}, are the sets of leaves of the two centers, of sizes $a$
and $b$. Hence Theorem~\ref{thm:pair} gives
$\gamma(F_2(S_{a,b}))\leq(a+b+2)-3+a=2a+b-1$, and by
Theorem~\ref{thm:doublestars} equality holds. We state this as follows.

\begin{corollary}\label{cor:sharp}
    The upper bound of Theorem~\ref{thm:pair} is attained by every double star
    $S_{a,b}$ with $2\leq a\leq b$; in particular, it cannot be improved in
    general.
\end{corollary}

%  %  %  %  %  %  %  %  %  %  %  %  %  %  %  %  %  %  %  %  %  %  %  %  %  %  % 
%                         S  U  B  S  E  C  T  I  O  N                        %
%  %  %  %  %  %  %  %  %  %  %  %  %  %  %  %  %  %  %  %  %  %  %  %  %  %  % 

\subsection{A lower bound, and the graphs with the maximum number of edges}

We now turn to the opposite extreme: graphs with domination number at least two
and the maximum possible number of edges.  The following lemma characterizes the
graphs attaining this maximum.
% We now turn to the opposite extreme: the graphs with domination number two and
% the maximum possible number of edges. These graphs are obtained from the
% complete graph $K_n$ by removing a maximum matching and for $n$ even, and,
% when $n$ is odd, one further edge. The following lemma determines them.

\begin{lemma}\label{lemma_maximum_edges}
    Let $G$ be a graph on $n$ vertices such that $\gamma(G)\geq 2$. Then
    \[
        |E(G)|
        \leq
        \binom{n}{2}-\left\lceil\frac n2\right\rceil,
    \]
    and equality holds if and only if $\overline{G}\cong\frac n2K_2$ when $n$ is
    even, and $\overline{G}\cong P_3\cup\frac{n-3}{2}K_2$ when $n$ is odd. In
    both cases $\gamma(G)=2$.
\end{lemma}

\begin{proof}
    Since $\gamma(G)\geq2$, $\Delta(G)\leq n-2$, and therefore every vertex of
    $G$ is incident with at least one edge of $\overline{G}$, which implies that
    $\overline{G}$ has no isolated vertices. A graph on $n$ vertices without
    isolated vertices has at most $\lfloor n/2\rfloor$ connected components,
    each of which has at least as many edges as its number of vertices minus
    one. Hence
    \[
        |E(\overline{G})|
        \geq n-\left\lfloor\frac n2\right\rfloor=\left\lceil\frac n2\right\rceil,
    \]
    and consequently
    \[
        |E(G)|
        =
        \binom n2-|E(\overline{G})|
        \leq
        \binom n2-\left\lceil\frac n2\right\rceil.
    \]
    Moreover, equality holds if and only if $\overline{G}$ has exactly $\lfloor
    n/2\rfloor$ components, all of them trees. Since each component has at least
    two vertices, this happens if and only if every component of $\overline{G}$
    is a $K_2$ when $n$ is even, and all but one are $K_2$ and the remaining one
    is a tree on three vertices, that is, a $P_3$, when $n$ is odd.

    It remains to see that in both cases $\gamma(G)=2$. Since $\overline{G}$ has
    no isolated vertices, $G$ has domination number at least two. On the other
    hand, one can choose two vertices belonging to a common edge of
    $\overline{G}$ so that every other vertex is adjacent in $G$ to at least one
    of them: for $n$ even any edge of $\overline{G}$ works, and for $n$ odd we
    take the middle vertex $x_2$ of the $P_3=x_1x_2x_3$ together with one of its
    two neighbors in $\overline{G}$, say $x_1$. Thus these two vertices form a
    dominating set of $G$, and therefore $\gamma(G)=2$.
\end{proof}

Since $F_k(H)$ is a spanning subgraph of $F_k(G)$ whenever $H$ is a spanning subgraph of $G$, we have the following.

\begin{remark}
    \label{rem:monotone}
    If $H$ is a spanning subgraph of a graph $G$, then $\gamma(F_k(G)) \le
    \gamma(F_k(H))$. In particular, if $G$ is a graph on $n$ vertices, then
    $\gamma(F_k(K_n)) \le \gamma(F_k(G))$.
    % Let $H$ be a spanning subgraph of $G$. Then $F_k(H)$ is a spanning
    % subgraph of $F_k(G)$, and therefore $\gamma(F_k(H))\geq\gamma(F_k(G))$. In
    % particular, $\gamma(F_k(G))\geq\gamma(F_k(K_n))$ for every graph $G$ on
    % $n$ vertices.
\end{remark}

Our next result gives a lower bound for $\gamma(F_2(G))$ which holds for every
graph without universal vertices, and shows that it is attained by the graphs of
Lemma~\ref{lemma_maximum_edges}.

\thmdense*

\begin{proof}
    We first prove the lower bound. In \cite{domination_token} it is shown that
    $\gamma(F_2(K_n)) =\lfloor n/2\rfloor$. Hence, by Remark~\ref{rem:monotone},
    \[
        \gamma(F_2(G))\geq\gamma(F_2(K_n))=\left\lfloor\frac{n}{2}\right\rfloor,
    \]
    which is the desired bound when $n$ is even. Suppose now that $n$ is odd,
    and suppose, aiming for a contradiction, that $\gamma(F_2(G))=\frac{n-1}2$,
    and let $D'$ be a dominating set of $F_2(G)$ with $|D'|=\frac{n-1}{2}$. By
    Remark~\ref{rem:monotone}, $D'$ is also a dominating set of $F_2(K_n)$. We
    now claim that $D'$ must be a matching of $K_n$ covering exactly $n-1$
    vertices.  Indeed, if two vertices of $K_n$ were not contained in any member
    of $D'$, then the corresponding $2$-set would be disjoint from every element
    of $D'$ and hence would not be dominated in $F_2(K_n)$. Since $D'$ contains
    exactly $(n-1)/2$ pairs, it follows that its members are pairwise disjoint
    and cover all but one vertex.

    Let $w$ be the unique vertex not covered by $D'$. Since $\gamma(G)\geq2$,
    the vertex $w$ is not universal in $G$.  Therefore, there exists a vertex
    $z$ such that $wz\notin E(G)$. Since $z$ is covered by $D'$, there exists a
    vertex $z'$ such that $\{z,z'\}\in D'$. Consider now the vertex
    $A=\{w,z'\}$. Since $w$ is not covered by $D'$, the only member of $D'$ that
    has nonempty intersection with $A$ is $\{z,z'\}$. However,
    $A\not\sim\{z,z'\}$, since adjacency between these two configurations would
    require $wz\in E(G)$. All other elements of $D'$ are disjoint from $A$, and
    hence cannot be adjacent to $A$ in $F_2(G)$. Therefore $A$ is not dominated
    by $D'$, a contradiction. Consequently
    \[
        \gamma(F_2(G))\geq\frac{n+1}{2}=\left\lceil\frac n2\right\rceil.
    \]

    We now show that equality holds for the graphs with the maximum number of
    edges. We consider the two cases separately.

    \medskip

    \noindent\textbf{Case 1: $n$ is even.} In this case, $G=K_n-M$, where $M$ is
    a perfect matching. Let
    \[
        M=\{\{u_1,v_1\},\ldots,\{u_{n/2},v_{n/2}\}\}
    \]
    and regard $D \coloneqq M$ as a subset of $V(F_2(G))$. We claim that $D$ is
    a dominating set of $F_2(G)$. Indeed, let $A=\{a,b\}\notin D$. Suppose,
    without loss of generality, that $a=u_i$. Since $A\neq\{u_i,v_i\}$, we have
    $b\neq v_i$. The only non-neighbor of $v_i$ in $G$ is $u_i$, and hence
    $bv_i\in E(G)$. Therefore, $\{a,b\}\sim\{u_i,v_i\}\in D$. Thus $D$ dominates
    $F_2(G)$, and consequently
    \[
        \gamma(F_2(G))\leq\frac{n}{2}.
    \]

    \medskip

    \noindent\textbf{Case 2: $n$ is odd.} Write
    \[
        \overline{G}
        =
        P_3\cup\frac{n-3}{2}K_2,
    \]
    where $P_3=x_1x_2x_3$, and the remaining edges of $\overline{G}$ are
    $u_1v_1,\ldots,u_{\frac{n-3}{2}}v_{\frac{n-3}{2}}$. Let
    \[
        D=
        \bigl\{\{x_1,x_2\},\{x_1,x_3\}\bigr\}
        \cup
        \left\{
        \{u_i,v_i\}:
        1\leq i\leq\frac{n-3}{2}
        \right\}.
        \]
        Thus
        \[
        |D|
        =
        2+\frac{n-3}{2}
        =
        \frac{n+1}{2}.
    \]

    We claim that $D$ is a dominating set of $F_2(G)$. First, suppose that a
    configuration $A$ contains some $u_i$ or $v_i$, say $A=\{u_i,w\}$. If
    $w=v_i$, then $A\in D$. Otherwise, $v_iw\in E(G)$, since the only
    non-neighbor of $v_i$ in $G$ is $u_i$. Hence $A\sim\{u_i,v_i\}\in D$. It
    remains to consider configurations contained in $\{x_1,x_2,x_3\}$. The
    configurations $\{x_1,x_2\}$ and $\{x_1,x_3\}$ belong to $D$. Also note that
    $\{x_2,x_3\}\sim\{x_1,x_2\}$, since $x_1x_3\in E(G)$. Thus $D$ is a
    dominating set, and
    \[
        \gamma(F_2(G))
        \leq
        \frac{n+1}{2}.
    \]
    Hence, in both cases, $\gamma(F_2(G))=\lceil n/2\rceil$.
\end{proof}

% % % % % % % % % % % % % % % % % % % % % % % % % % % % % % % % % % % % % % % %
%  %  %  %  %  %  %  %  %  %  %  %  %  %  %  %  %  %  %  %  %  %  %  %  %  %  % 
%                                                                             %
%                             S  E  C  T  I  O  N                             %
%                                                                             %
%  %  %  %  %  %  %  %  %  %  %  %  %  %  %  %  %  %  %  %  %  %  %  %  %  %  % 
% % % % % % % % % % % % % % % % % % % % % % % % % % % % % % % % % % % % % % % %

\section{Complete bipartite graphs}\label{sec:bipartite}

We now consider another family of graphs with domination number two, the
complete bipartite graphs $K_{r,t}$. In what follows, we assume that $2\leq
r\leq t$ are integers. Let $n=r+t$. Unless otherwise stated, we denote the
vertex bipartition of $K_{r,t}$ by
\[
    V_1 \coloneqq \{a_1,\dots,a_r\}\quad \text{and}\quad V_2 \coloneqq \{b_1,\dots,b_t\}.
\]

%  %  %  %  %  %  %  %  %  %  %  %  %  %  %  %  %  %  %  %  %  %  %  %  %  %  % 
%                         S  U  B  S  E  C  T  I  O  N                        %
%  %  %  %  %  %  %  %  %  %  %  %  %  %  %  %  %  %  %  %  %  %  %  %  %  %  % 

\subsection{Two tokens}\label{sec:bip2}

We begin by giving a general upper bound for the domination number of
$F_2(K_{r,t})$.

\begin{proposition}\label{prop_uppbnd}
    Let $2\leq r\leq t$. Then
    \[
        \gamma(F_2(K_{r,t}))
        \leq
        \begin{cases}
        t,
            & \text{if } r=2,3;\\[2mm]
        t+\dfrac{r-3}{2},
            & \text{if } r\geq 5 \text{ is odd and } t=r+1;\\[3mm]
        t+\left\lceil\dfrac{r}{2}\right\rceil-1,
            & \text{otherwise}.
        \end{cases}
    \]
\end{proposition}

\begin{proof}
    We first consider the cases $r=2$ and $r=3$.

    \medskip

    \begin{itemize}
        \item Suppose that $r=2$. Define
        \[
            D=
            \bigl\{\{a_1,a_2\}\bigr\}
            \cup
            \bigl\{\{a_1,b_i\}:1\leq i\leq t-1\bigr\}.
        \]
        Note that $|D|=1+(t-1)=t$. Every vertex of $F_2(K_{2,t})$ that contains one
        element of each part of the bipartition $V_1,V_2$ is adjacent to
        $\{a_1,a_2\}$.  Moreover, if $\{b_i,b_j\}\subseteq V_2$, then at least
        one of $i,j$ is different from $t$, say $i\neq t$, and hence is adjacent
        to $\{a_1,b_i\}\in D$. Thus $D$ is a dominating set and
        \[
            \gamma(F_2(K_{2,t}))\leq t.
        \]

        \item Suppose that $r=3$. Let
        \[
            D=
            \bigl\{\{a_1,b_i\}:1\leq i\leq t-2\bigr\}
            \cup
            \bigl\{\{a_2,a_3\},\{b_{t-1},b_t\}\bigr\}.
        \]
        Note again that $|D|=t$. If $A\in\binom{V_1}{2}$, then either
        $A=\{a_2,a_3\}\in D$, or $a_1\in A$, in which case $A\sim\{a_1,b_1\}\in
        D$.  If $A=\{b_i,b_j\}\in\binom{V_2}{2}$, then either
        $A=\{b_{t-1},b_t\}\in D$, or at least one of its elements, say $b_i$, is
        such that $i\leq t-2$, and therefore $A\sim\{a_1,b_i\}\in D$. Finally,
        consider a vertex $A=\{a_i,b_j\}$. If $i=1$ and $j\leq t-2$, then $A\in
        D$. If $i=1$ and $j\in\{t-1,t\}$, then $A\sim\{b_{t-1},b_t\}\in D$. And,
        if $i\in\{2,3\}$, then $A\sim\{a_2,a_3\}\in D$. Hence $D$ dominates
        $F_2(K_{3,t})$, and thus
        \[
            \gamma(F_2(K_{3,t}))\leq t.
        \]

        \item Suppose that $r\geq4$.  Let $\mathcal M\subseteq\binom{V_1}{2}$ be the following family of $\lceil r/2\rceil$ pairs, whose union is $V_1$:
        \[
            \mathcal M= \begin{cases}
            \bigl\{	\{a_1,a_2\},\{a_3,a_4\},\ldots,\{a_{r-1},a_r\}\bigr\} &
            \text{if $r$ is even,} \\
            \bigl\{ \{a_1,a_2\},\{a_3,a_4\},\ldots,
            \{a_{r-2},a_{r-1}\},\{a_{r-1},a_r\}\bigr\} &
            \text{if $r$ is odd.}
            \end{cases}
        \]

        In either case,
        \[
            |\mathcal M|=\left\lceil\frac r2\right\rceil.
        \]
        Let
        \[
            \mathcal C \coloneqq \bigl\{\{a_i,b_i\}:1\leq i\leq r-1\bigr\}
            \cup
            \bigl\{\{a_1,b_j\}:r\leq j\leq t-1\bigr\}
        \]
        (the second set is empty if $t=r$), and define the set
        \[
            D \coloneqq \mathcal M\cup\mathcal C.
        \]
        We claim that $D$ is a dominating set of $F_2(K_{r,t})$. Note that for
        any vertex $A=\{a_i,b_j\}$ of $F_2(K_{r,t})$ there is some $P\in\mathcal
        M$ containing $a_i$, and therefore $A\sim P$. Now, let
        $A=\{a_i,a_j\}\subseteq V_1$. In this case, at least one of $i,j$ is at
        most $r-1$, say $i\leq r-1$.  Thus, $A\sim\{a_i,b_i\}\in\mathcal C.$
        Finally, let $A=\{b_i,b_j\}\subseteq V_2$. Note that at least one of
        $i,j$ is at most $t-1$, say $i\leq t-1$.  Thus, if $i\leq r-1$, then
        $A\sim\{a_i,b_i\}\in\mathcal C,$ whereas if $r\leq i\leq t-1$, then
        $A\sim\{a_1,b_i\}\in\mathcal C.$ Therefore $D$ is a dominating set.
        Moreover,
        \[
            |\mathcal C|
            =(r-1)+(t-r)=t-1,
        \]
        and hence
        \[
            \gamma(F_2(K_{r,t}))\leq |D|=t+\left\lceil\frac r2\right\rceil-1.
        \]

        Suppose that $r\geq5$ is odd and $t=r+1$. In this case we can improve
        the last upper bound. Let
        \[
            \mathcal N=\bigl\{
            \{b_1,b_2\},\{b_3,b_4\},\ldots,
            \{b_{r-2},b_{r-1}\},\{b_r,b_{r+1}\}
            \bigr\}
        \]
        be a perfect matching of $V_2$ viewed as a subset of
        $V(F_2(K_{r,r+1}))$, and let
        \[
            \mathcal C=\bigl\{\{a_i,b_i\}:1\leq i\leq r-1\bigr\}.
        \]
        Finally, define
        \[
            D \coloneqq \mathcal N\cup\mathcal C.
        \]
        We claim that $D$ is a dominating set of $F_2(K_{r,r+1})$.  Indeed,
        every vertex $A=\{a_i,b_j\}$ is dominated by the unique pair of
        $\mathcal N$ containing $b_j$. If $\{a_i,a_j\}\subseteq V_1$, then at
        least one of its elements, say $a_i$, has $i\leq r-1$, and therefore
        $\{a_i,a_j\}\sim\{a_i,b_i\}\in\mathcal C.$ Finally, consider
        $\{b_i,b_j\}\subseteq V_2$. If one of its indices, say $i$, satisfies
        $i\leq r-1$, then
        \[
            \{b_i,b_j\}\sim\{a_i,b_i\}\in\mathcal C.
        \]
        The only pair of $V_2$ containing neither of $b_1,\ldots,b_{r-1}$ is
        $\{b_r,b_{r+1}\}$, which belongs to $\mathcal N$. Hence, $D$ is a
        dominating set and
        \[
            |D|=|\mathcal N|+|\mathcal C|=
            \frac{r+1}{2}+r-1=
            \frac{3r-1}{2}=
            (r+1)+\frac{r-3}{2}=
            t+\frac{r-3}{2}.
        \]
    \end{itemize}
    This completes the proof.
\end{proof}

The bounds of Proposition~\ref{prop_uppbnd} are in fact the exact values of
$\gamma(F_2(K_{r,t}))$: Acosta Troncoso, Caldera Caldera, and
Rivera~\cite{AcostaCalderaRivera2026} and Chinman, Cortijo, Pedroni, Seong, and
Stenger~\cite{Chinman2026} proved independently that
$\gamma(F_2(K_{r,t}))=\lceil r/2\rceil+t-2$ if $r=3$, or $r\geq5$ is odd and
$t=r+1$, and $\gamma(F_2(K_{r,t}))=\lceil r/2\rceil+t-1$ otherwise, for all
$2\leq r\leq t$. For $r=2$ we include our proof of the lower bound, which was obtained independently.

\begin{theorem}\label{thm:K2n}
    For every integer $n\geq 4$,
    \[
        \gamma(F_2(K_{2,n-2}))=n-2.
    \]
\end{theorem}

\begin{proof}
    The upper bound  is given by Proposition \ref{prop_uppbnd}, so it remains to
    show the lower bound.  To do this, suppose, to the contrary, that there
    exists a dominating set $D$ with at most $n-3$ vertices. Partition the
    vertices of $F_2(K_{2,n-2})$ into the sets
    \[
        U=\{\{a_1,a_2\}\},\quad C=
        \{\{a_i,b_j\}:i\in[2], j\in [n-2]\}, \quad  W=
        \{\{b_i,b_j\}:1\leq i<j\leq n-2\}.
    \]

    We distinguish two cases.

    \medskip

    \noindent
    \textbf{Case 1.}
    Suppose that $\{a_1,a_2\}\in D$. Let $x=|D\cap C|$. Then
    \[
        |D\cap W|
        \leq
        n-4-x.
    \]

    Each vertex of $C$ dominates exactly $n-3$ vertices of $W$. Moreover, the
    number of vertices of $W$ dominated by $x$ vertices of $C$ is maximized when
    the corresponding indices of $V_2$ are distinct, in which case this number
    is $x(n-3)-\binom{x}{2}$.  Therefore, at most
    \[
        f(x)
        =
        x(n-3)
        -
        \binom{x}{2}
        +
        n-4-x
    \]
    vertices of $W$ are either contained in $D$ or adjacent to a vertex of
    $D\cap C$. Since
    \[
        f(x+1)-f(x)
        =
        n-4-x\geq0,
    \]
    the maximum of $f$ is attained at $x=n-4$. Consequently,
    \[
        f(n-4)
        =
        \binom{n-2}{2}-1
        <
        |W|.
    \]
    Hence some vertex of $W$ is not dominated, a contradiction.

    \medskip

    \noindent
    \textbf{Case 2.}
    Suppose that $\{a_1,a_2\}\notin D$. Let $s$ be the number of indices
    $j\in\{1,\ldots,n-2\}$ for which at least one of the vertices $\{a_1,b_j\}$
    or $\{a_2,b_j\}$, belongs to $D$, and let $d$ be the number of indices for
    which both belong to $D$. Thus
    \[|D\cap C|=s+d.  \] Let
    \[
        R=
        \{b_j:
        \{a_1,b_j\},
        \{a_2,b_j\}
        \notin D\}.
    \]
    Then
    \[
        |R|=n-2-s.
    \]
    If $|R|\geq2$, every pair contained in $R$ must belong to $D$, and therefore
    \[
        |D\cap W|
        \geq
        \binom{n-2-s}{2}.
    \]

    Furthermore, there are $s-d$ indices represented exactly once in $D\cap C$.
    Each of them must be incident with a vertex of $D\cap W$, and every vertex
    of $W$ contains at most two such indices.  Hence
    \[
        |D\cap W|
        \geq
        \binom{n-2-s}{2}
        +
        \left\lceil
            \frac{s-d}{2}
        \right\rceil.
    \]

    \noindent Consequently,
    \[
        |D|
        \geq
        F_n(s,d), \qquad \text{where}\qquad
        F_n(s,d)
         \coloneqq 
        s+d+
        \binom{n-2-s}{2}
        +
        \left\lceil\frac{s-d}{2}\right\rceil.
    \]

    We claim that $F_n(s,d)\geq n-2$. First observe that, for fixed $s$,
    \[
        F_n(s,d+1)-F_n(s,d)=
        1+
        \left\lceil\frac{s-d-1}{2}\right\rceil-  \left\lceil\frac{s-d}{2}\right\rceil
        \in\{0,1\}.
    \]
    Thus, $F_n(s,d)$ is nondecreasing as a function of $d$, and hence
    \[
        F_n(s,d)\geq F_n(s,0).
    \]

    Set $r=n-2-s$. Since $\{a_1,a_2\}\notin D$, at least one vertex of $C$ must
    belong to $D$ in order to dominate $\{a_1,a_2\}$. Hence $s\geq1$. For
    $r\geq2$, we have
    \[
        F_n(s,0)=
        n-2-r+\binom{r}{2}
        +
        \left\lceil\frac{n-2-r}{2}\right\rceil=
        n-2+
        \frac{r(r-3)}{2}
        +
        \left\lceil\frac{n-2-r}{2}\right\rceil.
    \]
    If $r\geq3$, then both terms after $n-2$ are nonnegative, and therefore
    \[
        F_n(s,0)\geq n-2.
    \]
    If $r=2$, then $s=n-4\geq1$ (for $n=4$ this case does not arise, as
    $s\geq1$), and therefore
    \[
        F_n(s,0)=
        n-2-1+
        \left\lceil\frac{n-4}{2}\right\rceil\geq n-2.
    \]
    Thus,
    \[
        |D|\geq F_n(s,d)\geq n-2,
    \]
    contradicting the assumption that $|D|\leq n-3$.

    Now we analyze the remaining cases $|R|=0$ and $|R|=1$, for which we also
    have that $|D|\geq n-2$.  Indeed, if $|R|=0$, then
    \[
    s=n-2\qquad \text{and}\qquad|D|\geq|D\cap C|=s+d\geq n-2.
    \]
    Now, if $|R|=1$, say $R = \{b_{j_0} \}$, then the configuration
    $\{a_1,b_{j_0}\}$ requires some $\{b_k,b_{j_0}\}\in D$, so that there are
    $s-d+1$ indices that require a pair in $W$, and since each vertex of $W$
    contains at most two of them, we have
    \[
        |D|\ge(n-3)+d+\Bigl\lceil\tfrac{n-2-d}{2}\Bigr\rceil\ge n-2.
    \]  

    Hence
    \[
        \gamma(F_2(K_{2,n-2}))
        \geq
        n-2,
    \]
    as desired.
\end{proof}

%  %  %  %  %  %  %  %  %  %  %  %  %  %  %  %  %  %  %  %  %  %  %  %  %  %  % 
%                         S  U  B  S  E  C  T  I  O  N                        %
%  %  %  %  %  %  %  %  %  %  %  %  %  %  %  %  %  %  %  %  %  %  %  %  %  %  % 

\subsection{Three or more tokens}\label{sec:bip3}

We now consider $k\geq3$. We first recall some facts on the Tur\'an numbers
$T(m,k,k-1)$ defined in the introduction. The ratio $T(m,k,k-1)/\binom m{k-1}$ is
non-decreasing in $m$, by a standard averaging argument over all $m$-subsets of a
larger ground set, so the limit $t_k$ exists. For $k=3$ a Tur\'an $(m,3,2)$-family is the edge set of a graph whose
complement is triangle-free, so by Mantel's theorem~\cite{Mantel1907}
\begin{equation}\label{eq:mantel}
    T(m,3,2) = \binom m2-\left\lfloor\frac{m^2}{4}\right\rfloor
    =\binom{\lceil m/2\rceil}{2}+\binom{\lfloor m/2\rfloor}{2},
\end{equation}
attained precisely when the family consists of all pairs inside each of the two
classes of a balanced bipartition of $[m]$; thus $t_3=\frac12$.
Bounds on $t_k$ for $k\geq4$ are discussed after the proof of Theorem~\ref{thm:turan}.

We recall the following well-known bound, already mentioned in the introduction.

\begin{proposition}[\cite{HHS98}]\label{prop_gen_lwb}
    Let $G$ be a graph on $n$ vertices. Then,
    \[
        \gamma(G)\ge \frac{n}{\Delta(G)+1}.
    \]
\end{proposition}

\noindent We now give a general lower bound for $\gamma(F_k(K_{r,t}))$ based on
Proposition~\ref{prop_gen_lwb}, for which we compute the maximum degree of
$F_k(K_{r,t})$.

\begin{proposition} \label{prop_lwbnd_token}
If $2\leq r\leq t$ and $2\leq k\leq n/2$, then
\[
\gamma(F_k(K_{r,t}))
\geq
\frac{\binom{n}{k}}{\Delta(F_k(K_{r,t}))+1},
\]
where
\[
\Delta(F_k(K_{r,t}))
=
\max_{0\leq i\leq \min\{k,r\}}
\bigl\{i(t-k+i)+(k-i)(r-i)\bigr\}.
\]
% In particular, if $k\leq r$, then
% \[
% \gamma(F_k(K_{r,t}))
% \geq
% \frac{\binom{n}{k}}{kt+1}.
% \]
Moreover, if $3\leq k\leq r\leq t$, then
\[
\gamma(F_k(K_{r,t}))
\geq
\max\left\{
\left\lceil\frac{\binom{r+t}{k}}{kt+1}\right\rceil,
T(r,k,k-1)+T(t,k,k-1)
\right\}.
\]
\end{proposition}

\begin{proof}
We first establish the bound in terms of the
maximum degree. Let $A\in V(F_k(K_{r,t}))$ and
suppose that
\[
|A\cap V_1|=i,
\qquad
|A\cap V_2|=k-i.
\]
Each of the $i$ tokens on $V_1$ can be moved to
any of the $t-(k-i)$ unoccupied vertices of $V_2$.
Similarly, each of the $k-i$ tokens on $V_2$ can
be moved to any of the $r-i$ unoccupied vertices
of $V_1$. Hence
\[
\deg(A)
=
i(t-k+i)+(k-i)(r-i).
\]
Since $k\leq n/2\leq t$, the possible values of
$i$ are $0\leq i\leq\min\{k,r\}$. Therefore,
\[
\Delta(F_k(K_{r,t}))
=
\max_{0\leq i\leq\min\{k,r\}}
\bigl\{i(t-k+i)+(k-i)(r-i)\bigr\}.
\]
The first lower bound now follows from
Proposition~\ref{prop_gen_lwb}.

Now, if $2\leq k\leq r\leq t$, then all values $0\leq i\leq k$
are feasible. Moreover,
\[
i(t-k+i)+(k-i)(r-i)=
2i^2+(t-r-2k)i+kr
\]
is a convex quadratic function of $i$, so its
maximum on $[0,k]$ is attained at an endpoint.
Since the values at $i=0$ and $i=k$ are $kr$
and $kt$, respectively, and $r\leq t$, we obtain
\[
\Delta(F_k(K_{r,t}))=kt.
\]
Consequently,
\[
\gamma(F_k(K_{r,t}))
\geq
\left\lceil
\frac{\binom{r+t}{k}}{kt+1}
\right\rceil.
\]

% Now suppose that $3\leq k\leq r\leq t$.
Now, let $D$ be a dominating set of $F_k(K_{r,t})$ and assume that $3\leq k\leq r\leq t$. Consider the
configurations contained entirely in $V_1$. Such a configuration
$A\in\binom{V_1}{k}$ can only be dominated by itself or by a configuration
containing exactly $k-1$ vertices of $V_1$ and one vertex of $V_2$. Define
\[
D_1=
\left\{
A\in D:|A\cap V_1|\geq k-1
\right\}.
\]
For each $A\in D_1$, choose a
$(k-1)$-subset $P_A\subseteq A\cap V_1$,
and let
\[
\mathcal{F}_1=\{P_A:A\in D_1\}.
\]
Since $D$ dominates every configuration
in $\binom{V_1}{k}$, each such configuration
contains some member of $\mathcal{F}_1$.
Thus $\mathcal{F}_1$ is a Tur\'an
$(r,k,k-1)$-family, and hence
\[
|D_1|\geq|\mathcal{F}_1|
\geq T(r,k,k-1).
\]
Analogously, considering configurations
contained entirely in $V_2$, define
\[
D_2=
\left\{
A\in D:|A\cap V_2|\geq k-1
\right\}.
\]
The same argument yields
\[
|D_2|\geq T(t,k,k-1).
\]

Since $k\geq3$, no configuration can contain
at least $k-1$ vertices from both $V_1$
and $V_2$. Therefore $D_1\cap D_2=\varnothing$,
and consequently
\[
|D|\geq |D_1|+|D_2|
\geq T(r,k,k-1)+T(t,k,k-1).
\]

Combining this inequality with the
maximum-degree bound gives
\[
\gamma(F_k(K_{r,t}))
\geq
\max\left\{
\left\lceil
\frac{\binom{r+t}{k}}{kt+1}
\right\rceil,
T(r,k,k-1)+T(t,k,k-1)
\right\},
\]
as desired.
\end{proof}

\begin{remark}
The two lower bounds in Proposition~\ref{prop_lwbnd_token} are complementary:
neither one is always stronger than the other. For instance, when $k=3$ and
$r=t=4$, the maximum-degree bound gives $\gamma(F_3(K_{4,4}))\geq 5$, whereas
the Tur\'an bound gives only $4$; conversely, for $k=r=3$ and $t=10$ the maximum-degree bound gives $10$, whereas the Tur\'an bound gives $T(3,3,2)+T(10,3,2)=21$.
\end{remark}

We now consider the graphs $K_{2,n-2}$. Let $V_1=\{x_1,x_2\}$ and $V_2=[n-2]$ be the
vertex bipartition of $K_{2,n-2}$; throughout, $x$ denotes a vertex of $V_1$ and
$x'$ the other vertex of $V_1$, while $y,y'$ denote vertices of $V_2$.  For
$i\in\{0,1,2\}$ we say that a configuration $A$ of $F_k(K_{2,n-2})$ is a
\emph{type-$i$ configuration} if $|A\cap V_1|=i$.  Since every edge of
$K_{2,n-2}$ joins $V_1$ and $V_2$, moving a token changes the type by one.  In
particular, two configurations of the same type are never adjacent. 
We now describe the neighbors of each configuration; here and below, $P$ always denotes a $(k-1)$-subset and $Q$ a
$(k-2)$-subset of $V_2$.
\begin{itemize}
    \item The neighbors of a type-$0$ configuration $A$ are the type-$1$
    configurations $\{x\}\cup P$ with $P\subseteq A$.

    \item The neighbors of a type-$1$ configuration $\{x\}\cup P$ are the
    type-$0$ configurations $P\cup\{y\}$ with $y\in V_2\setminus P$, and the
    type-$2$ configurations $\{x_1,x_2\}\cup(P\setminus\{y\})$ with $y\in P$.

    \item The neighbors of a type-$2$ configuration $\{x_1,x_2\}\cup Q$ are the
    type-$1$ configurations $\{x\}\cup P$ with $Q\subseteq P$.
\end{itemize}

We need two more extremal notions. A family
$\mathcal{C}\subseteq\binom{[m]}{k-1}$ is a \emph{covering of the $(k-2)$-sets}
if every $(k-2)$-subset of $[m]$ is contained in a member of $\mathcal{C}$; the
minimum size of such a family is the \emph{covering number} $C(m,k-1,k-2)$.
Since every $(k-1)$-set contains exactly $k-1$ subsets of size $k-2$, we have
$C(m,k-1,k-2)\geq\frac{1}{k-1}\binom{m}{k-2}$. For fixed $k$, a theorem of
R\"odl~\cite{Rodl1985}, applied with blocks of size $k-1$ covering the
$(k-2)$-sets, shows that this bound is asymptotically tight:
\[
\lim_{m\to\infty}\frac{(k-1)\,C(m,k-1,k-2)}{\binom{m}{k-2}}=1,
\]
that is,
\begin{equation}\label{eq:rodl}
C(m,k-1,k-2)=\bigl(1+o(1)\bigr)\frac{1}{k-1}\binom{m}{k-2}.
\end{equation}

For $k=4$, a covering of the $2$-sets is a family of triples containing every
pair, also called an \emph{$(m,3,2)$-covering}, and $C(m,3,2)=\left\lceil\frac
m3\left\lceil\frac{m-1}2\right\rceil\right\rceil$ for every
$m\geq3$~\cite{FortHedlund1958}; in particular,
\begin{equation}\label{eq:FH}
C(m,3,2)=\frac{m^2}{6}+O(m).
\end{equation}
Note also that a Tur\'an $(m,k-1,k-2)$-family is a family of $(k-2)$-sets such
that every $(k-1)$-set contains one of them; for $k=3$ it is a set of vertices
meeting every pair, so $T(m,2,1)=m-1$.

For $k=4$ we will use the following construction of Tur\'an~\cite{Turan1941}.
Let $m\geq6$, let $[m]=W_0\sqcup W_1\sqcup W_2$ be a partition into three
classes whose sizes differ by at most one, with indices taken modulo $3$, and
let $\mathcal{T}_m\subseteq\binom{[m]}{3}$ be the family of the triples that
either are contained in some class $W_i$, or have two vertices in $W_i$ and one
vertex in $W_{i-1}$ for some $i$. The complement
$\binom{[m]}{3}\setminus\mathcal{T}_m$ consists of the triples meeting the three
classes and of the triples with two vertices in $W_i$ and one in $W_{i+1}$; it
contains no copy of $K_4^{(3)}$ (see~\cite{Keevash2011}), so $\mathcal{T}_m$ is
a Tur\'an $(m,4,3)$-family, and $|\mathcal{T}_m|=(\frac49+o(1))\binom m3$. We
call $\mathcal{T}_m$ \emph{Tur\'an's construction}. Tur\'an's conjecture, in its
exact form, states that $|\mathcal{T}_m|=T(m,4,3)$. Moreover, every pair of
$[m]$ lies in a member of $\mathcal{T}_m$: a pair contained in $W_i$, or a pair
$\{u,w\}$ with $u\in W_i$ and $w\in W_{i-1}$, is contained in a triple with two
vertices in $W_i$ and one in $W_{i-1}$ (recall that $|W_i|\geq2$), and every
pair meeting two classes is of the second form for a suitable $i$.

We now prove Theorem~\ref{thm:turan}, which we restate for easy reference.

\thmturan*

\begin{proof}
Write $m\coloneqq n-2$, so that $V_2=[m]$ and $m\geq2k-2\geq k$.

We first prove the lower bound. Let $D$ be a dominating set of $F_k(K_{2,m})$,
and let $\mathcal{F}_{k-1}\subseteq\binom{V_2}{k-1}$ be the family consisting of
all $P$ such that $\{x_1\}\cup P\in D$ or $\{x_2\}\cup P\in D$, together with,
for each type-$0$ configuration $A\in D$, one arbitrarily chosen $(k-1)$-subset
$P_A$ of $A$. Every $k$-subset $A$ of $V_2$ contains a member of
$\mathcal{F}_{k-1}$: if $A\in D$, then $P_A\in\mathcal{F}_{k-1}$, and otherwise,
by the first item above, $A\sim\{x\}\cup P\in D$ for some $x\in V_1$ and some
$P\subseteq A$, so that $P\in\mathcal{F}_{k-1}$. That is, $\mathcal{F}_{k-1}$ is
a Tur\'an $(m,k,k-1)$-family. Since every member of $\mathcal{F}_{k-1}$ arises
from at least one element of $D$, we have $|D|\geq|\mathcal{F}_{k-1}|\geq
T(m,k,k-1)$.

For the upper bound, let $\mathcal{F}_{k-1}\subseteq\binom{V_2}{k-1}$ be a
Tur\'an $(m,k,k-1)$-family, let $\mathcal{F}_{k-2}\subseteq\binom{V_2}{k-2}$ be
a Tur\'an $(m,k-1,k-2)$-family, and let $\mathcal{C}\subseteq\binom{V_2}{k-1}$
be such that every $(k-2)$-subset of $V_2$ not in $\mathcal{F}_{k-2}$ is
contained in a member of $\mathcal{F}_{k-1}\cup\mathcal{C}$. Define
\[
D\coloneqq \bigl\{\{x_1\}\cup P:P\in\mathcal{F}_{k-1}\cup\mathcal{C}\bigr\}\ \cup\
\bigl\{\{x_1,x_2\}\cup Q:Q\in\mathcal{F}_{k-2}\bigr\}.
\]
Let $A$ be a configuration of $F_k(K_{2,m})$. If $A$ is a type-$0$
configuration, then $A$ contains some $P\in\mathcal{F}_{k-1}$, and
$A\sim\{x_1\}\cup P$, since their symmetric difference is $\{x_1,y\}$, where $y$
is the unique element of $A\setminus P$. If $A=\{x\}\cup P$ is a type-$1$
configuration, then $P$ contains some $Q\in\mathcal{F}_{k-2}$, and
$A\sim\{x_1,x_2\}\cup Q$, since their symmetric difference is $\{x',y\}$, where
$y$ is the unique element of $P\setminus Q$. Finally, let $A=\{x_1,x_2\}\cup Q$
be a type-$2$ configuration. If $Q\in\mathcal{F}_{k-2}$, then $A\in D$;
otherwise $Q\subseteq P$ for some $P\in\mathcal{F}_{k-1}\cup\mathcal{C}$, and
$A\sim\{x_1\}\cup P$, since their symmetric difference is $\{x_2,y\}$, where $y$
is the unique element of $P\setminus Q$. Hence $D$ is a dominating set, and
$|D|\leq|\mathcal{F}_{k-1}|+|\mathcal{F}_{k-2}|+|\mathcal{C}|$. Taking
$\mathcal{F}_{k-1}$ and $\mathcal{F}_{k-2}$ of minimum size and $\mathcal{C}$ a
covering of all the $(k-2)$-sets of size $C(m,k-1,k-2)$ gives the first upper
bound, and taking $\mathcal{C}=\varnothing$ gives the second
one.

For the last assertion, note that 
$T(m,k-1,k-2)=(t_{k-1}+o(1))\binom{m}{k-2}$,
$C(m,k-1,k-2)=(1+o(1))\frac1{k-1}\binom m{k-2}$ by~\eqref{eq:rodl}, and
$T(m,k,k-1)=(t_k+o(1))\binom{m}{k-1}=(t_k+o(1))\binom{n}{k-1}$.
\end{proof}

Taking complements, a family $\mathcal{F}_{k-1}\subseteq\binom{[m]}{k-1}$ is a
Tur\'an $(m,k,k-1)$-family if and only if
$\binom{[m]}{k-1}\setminus\mathcal{F}_{k-1}$ contains no copy of the complete
$(k-1)$-uniform hypergraph $K_k^{(k-1)}$ on $k$ vertices, that is, no $k$-set
all of whose $(k-1)$-subsets lie outside $\mathcal{F}_{k-1}$. Hence
\[
T(m,k,k-1)=\binom{m}{k-1}-\operatorname{ex}\bigl(m,K_k^{(k-1)}\bigr)
\qquad\text{and}\qquad
t_k=1-\pi\bigl(K_k^{(k-1)}\bigr),
\]
where $\operatorname{ex}$ and $\pi$ denote the Tur\'an number and the Tur\'an
density, so that $\gamma(F_k(K_{2,n-2}))$ is asymptotically determined by the
classical hypergraph Tur\'an problem for $K_k^{(k-1)}$. In particular, $t_k$ and $\pi(K_k^{(k-1)})$ coincide
only for $k=3$, where both are $\frac12$ by Mantel's theorem~\cite{Mantel1907}.
By results of de~Caen~\cite{deCaen1983} and Gao et al.~\cite{GaoEtAl2026}, the
latter improving a bound of Pikhurko~\cite{Pikhurko2024},
$\frac1{k-1}\leq t_k\leq\frac{4}{k+1}$, and $t_k$ is known exactly only for
$k=3$.  For $k=4$, Tur\'an's construction $\mathcal{T}_m$ gives
$t_4\leq\frac49$, and Tur\'an~\cite{Turan1941} conjectured that equality holds;
on the other hand, the recent bound $\pi(K_4^{(3)})\leq0.5578$ of Jeong et
al.~\cite{JeongEtAl2026} gives $t_4\geq0.4422$.  Consequently,
\[
    \bigl(0.4422-o(1)\bigr)\binom{n}{3}\ \leq\ \gamma\bigl(F_4(K_{2,n-2})\bigr)\ \leq\
    \Bigl(\frac49+o(1)\Bigr)\binom{n}{3},
\]
so that $\gamma(F_4(K_{2,n-2}))/\binom n3$ converges to
$t_4\in[0.4422,\frac49]$.  

The second-order term of the upper bound of Theorem~\ref{thm:turan} depends on the Tur\'an
and covering problems of the previous order; for $k=4$ these values are explicit,
by~\eqref{eq:FH} and~\eqref{eq:mantel}:
\begin{align*}
\gamma\bigl(F_4(K_{2,n-2})\bigr)\ &\leq\ T(n-2,4,3)+\binom{n-2}{2}-\left\lfloor\frac{(n-2)^2}{4}\right\rfloor
+\left\lceil\frac{n-2}{3}\left\lceil\frac{n-3}{2}\right\rceil\right\rceil\\
&=\ T(n-2,4,3)+\frac{5}{12}n^2+O(n).
\end{align*}

Recall that the second bound of Theorem~\ref{thm:turan} requires a Tur\'an
family $\mathcal{F}_{k-1}$ of $(k-1)$-sets and a Tur\'an family
$\mathcal{F}_{k-2}$ of $(k-2)$-sets such that every $(k-2)$-set not in
$\mathcal{F}_{k-2}$ lies in a member of $\mathcal{F}_{k-1}$. For $k=3$ and
$m\geq4$, this holds for the family of~\eqref{eq:mantel}, in which every vertex
lies in a pair. For $k=4$ and $m\geq6$, it holds for Tur\'an's construction
$\mathcal{T}_m$ and the family of~\eqref{eq:mantel}, since every pair of $[m]$
lies in a member of $\mathcal{T}_m$. Hence
$\gamma(F_4(K_{2,m}))\leq|\mathcal{T}_m|+T(m,3,2)$, which equals
$T(m,4,3)+T(m,3,2)$ if Tur\'an's conjecture holds.

For small $m$ this bound is close to the exact value. Since
$|\mathcal{T}_m|=T(m,4,3)$ for $6\leq m\leq9$, it gives
$\gamma(F_4(K_{2,m}))\leq12,21,32,46$ for these values of $m$, while a computer
search gives the exact values  $\gamma(F_4(K_{2,m}))=12,19,30$ for $m=6,7,8$.

For $k=3$ we can determine $\gamma(F_3(K_{2,n-2}))$ exactly. Here
$T(n-2,2,1)=n-3$, since a set of vertices meeting every pair of $[n-2]$ contains
all the vertices but at most one, and the hypothesis of the second bound of
Theorem~\ref{thm:turan} holds for the family of~\eqref{eq:mantel} when $n\geq6$;
hence $\gamma(F_3(K_{2,n-2}))\leq T(n-2,3,2)+n-3$. In the proof below we save
one more element, by replacing two type-$2$ configurations by a single type-$0$
configuration, and a supersaturation argument shows that the resulting bound is
the exact value.  We use the following two classical results on the number of
triangles of a graph with slightly more than $\lfloor n^2/4\rfloor$ edges.

\begin{lemma}\label{lem:supersat}
Let $G$ be a graph on $n$ vertices, $m$ edges, and $\operatorname{tr}(G)$
triangles. Then,
\begin{enumerate}[label=(\roman*)]
    \item (Nordhaus and Stewart~\cite{nordhaus1963triangles}; see also Moon and
    Moser~\cite{MoonMoser1962})
    $\operatorname{tr}(G)\geq\dfrac{m\,(4m-n^2)}{3n}$.

    \item (Lov\'asz and Simonovits~\cite{LovaszSimonovits1976}) If $m=\lfloor
    n^2/4\rfloor+s$ with $1\leq s<n/2$, then $\operatorname{tr}(G)\geq s\lfloor
    n/2\rfloor$.
\end{enumerate}
\end{lemma}

\begin{theorem}
    \label{thm:K2n3}
    For every $n\geq6$,
    \[
        \gamma\bigl(F_3(K_{2,n-2})\bigr) =
            \binom{n-2}{2}-\left\lfloor\frac{(n-2)^2}{4}\right\rfloor+n-4 = 
            T(n-2,3,2)+n-4 .
    \]
    In particular, $\gamma(F_3(K_{2,n-2}))=\frac{n^2}{4}+O(n)$.
\end{theorem}

\begin{proof}
    Write $m\coloneqq n-2\geq4$, and recall from~\eqref{eq:mantel} that
    $T(m,3,2)=\binom m2-\lfloor m^2/4\rfloor$.

    \emph{Upper bound.} Let $B_1\sqcup B_2$ be a balanced bipartition of $V_2$
    and let $\mathcal{F}_2\coloneqq \binom{B_1}{2}\cup\binom{B_2}{2}$, so that
    $|\mathcal{F}_2|=T(m,3,2)$ and, as $m\geq4$, every vertex of $V_2$ lies in a
    pair of $\mathcal{F}_2$. Let $T_0\coloneqq \{1,2,3\}\subseteq V_2$, and
    define
    \[
    D\coloneqq \bigl\{\{x_1\}\cup P:P\in\mathcal{F}_2\bigr\}\ \cup\ \bigl\{\{x_1,x_2,y\}:y\in V_2\setminus
    T_0\bigr\}\ \cup\ \{T_0\}.
    \]
    Every type-$0$ configuration $A$ has two vertices in the same class $B_j$,
    so it contains a pair $P\in\mathcal{F}_2$, and $A\sim\{x_1\}\cup P$.  A
    type-$2$ configuration $\{x_1,x_2,y\}$ belongs to $D$ if $y\notin T_0$; if
    $y\in T_0$, then $y$ lies in a pair $P\in\mathcal{F}_2$, and
    $\{x_1,x_2,y\}\sim\{x_1\}\cup P$. Finally, let $\{x,y,y'\}$ be a type-$1$
    configuration. If one of $y,y'$, say $y'$, is not in $T_0$, then
    $\{x,y,y'\}\sim\{x_1,x_2,y'\}\in D$; otherwise $\{y,y'\}\subseteq T_0$ and
    $\{x,y,y'\}\sim T_0$.  Hence $D$ is a dominating set, and
    $|D|=T(m,3,2)+(m-3)+1=T(n-2,3,2)+n-4$.

    \emph{Lower bound.} For $4\leq m\leq9$, that is, for $6\leq n\leq11$, the
    lower bound was verified by computer; the values of $\gamma(F_3(K_{2,m}))$
    for $m=4,\dots,9$ are $4,7,10,14,18,23$, which are those given by the
    formula. Hence we may assume that $m\geq10$.

    Let $D$ be a dominating set of $F_3(K_{2,m})$, and for $i\in\{0,1,2\}$ let
    $d_i$ be the number of type-$i$ configurations in $D$. We split $V_2=U\sqcup
    W$, where
    \[
    U\coloneqq \{y\in V_2:\{x_1,x_2,y\}\in D\}\qquad\text{and}\qquad W\coloneqq V_2\setminus U,
    \]
    and we write $w\coloneqq |W|$, so that $d_2=m-w$. As in the proof of
    Theorem~\ref{thm:turan}, let $\mathcal{F}_2$ be the set of pairs $P\subseteq
    V_2$ such that $\{x_1\}\cup P\in D$ or $\{x_2\}\cup P\in D$, and let
    $f\coloneqq |\mathcal{F}_2|$. Finally, let $\mathcal{R}$ be the set of pairs
    $P\subseteq W$ that are contained in no type-$0$ element of $D$, and let
    $r\coloneqq |\mathcal{R}|$.  The proof has three steps. In the first one we
    bound $|D|$ from below in terms of $f$, $w$ and the quantity $d_0+r$; in the
    second one we give two lower bounds for $d_0+r$; and in the third one we
    combine them, using Lemma~\ref{lem:supersat} to bound $d_0$ when
    $\mathcal{F}_2$ is smaller than a Tur\'an family of minimum size.

    \emph{Step 1.} Let $P\in\mathcal{R}$. Since $P\subseteq W$, no type-$2$
    element of $D$ dominates $\{x_1\}\cup P$ or $\{x_2\}\cup P$, and by the
    definition of $\mathcal{R}$ no type-$0$ element of $D$ does either; hence
    both $\{x_1\}\cup P$ and $\{x_2\}\cup P$ belong to $D$. In particular,
    $\mathcal{R}\subseteq\mathcal{F}_2$. Thus each pair of
    $\mathcal{F}_2\setminus\mathcal{R}$ accounts for at least one type-$1$
    element of $D$, and each pair of $\mathcal{R}$ for two, so that
    $d_1\geq(f-r)+2r=f+r$. Therefore
    \begin{equation}\label{eq:total}
    |D|\ =\ d_0+d_1+d_2\ \geq\ f+(d_0+r)+(m-w).
    \end{equation}

    \emph{Step 2.} First, let $H$ be the graph on $V_2$ whose edges are the
    pairs not in $\mathcal{F}_2$. A triangle of $H$ is a type-$0$ configuration
    containing no pair of $\mathcal{F}_2$, so it is dominated only by itself;
    therefore
    \begin{equation}\label{eq:d0}
    d_0\ \geq\ \operatorname{tr}(H).
    \end{equation}
    Second, each type-$0$ element of $D$ contains at most three pairs of $W$, so
    $r\geq\binom w2-3d_0$, and hence
    \[
    d_0+r\ \geq\ \max\Bigl\{d_0,\ \binom w2-2d_0\Bigr\}\ \geq\ \frac13\binom w2;
    \]
    indeed, if $d_0\geq\frac13\binom w2$ there is nothing to prove, and
    otherwise $\binom w2-2d_0>\binom w2-\frac23\binom w2=\frac13\binom w2$.
    Moreover, $\frac13\binom w2-(w-2)=\frac{(w-3)(w-4)}{6}\geq0$ for every
    integer $w$, since $(w-3)(w-4)$ is a product of two consecutive integers.
    Hence
    \begin{equation}\label{eq:w}
    d_0+r\ \geq\ \frac13\binom w2\ \geq\ w-2 .
    \end{equation}

    \emph{Step 3.} Let $s\coloneqq T(m,3,2)-f$, which measures how far
    $\mathcal{F}_2$ is from being a Tur\'an family of minimum size.  If
    $s\leq0$, then \eqref{eq:total} and \eqref{eq:w} give $|D|\geq
    f+(w-2)+(m-w)\geq T(m,3,2)+m-2$, as desired. So assume that $s\geq1$. Recall
    that $H$ is the graph on $V_2$ whose edges are the pairs not in
    $\mathcal{F}_2$.  Then $H$ has $\binom m2-f=\lfloor m^2/4\rfloor+s$ edges,
    and by~\eqref{eq:total} it suffices to prove that
    \begin{equation}\label{eq:goal}
    d_0+r\ \geq\ s+w-2 ,
    \end{equation}
    since then $|D|\geq f+s+(w-2)+(m-w)=T(m,3,2)+m-2$. Let $q\coloneqq \lfloor
    m/2\rfloor$, so that $q\geq5$.

    \emph{Case $s<m/2$.} By Lemma~\ref{lem:supersat}(ii) and~\eqref{eq:d0},
    $d_0\geq sq$. Since $d_0+r\geq d_0$, in the first two subcases below it is
    enough to show that $d_0\geq s+w-2$.
    \begin{itemize}
        \item If $w\leq q+1$, then, as $s\geq1$, $d_0\geq s+s(q-1)\geq s+q-1\geq
        s+w-2$.
        \item If $w\geq q+2$ and $s(q-1)\geq w-2$, then $d_0\geq s+s(q-1)\geq
        s+w-2$.
        \item If $w\geq q+2$ and $s(q-1)<w-2$, then $s<\frac{w-2}{q-1}$, so
        $s+w-2<(w-2)\bigl(1+\frac{1}{q-1}\bigr)=(w-2)\frac{q}{q-1}$, and
        by~\eqref{eq:w} it suffices to show that $\frac13\binom
        w2\geq(w-2)\frac{q}{q-1}$, that is, $h(w)\coloneqq
        \frac{w(w-1)}{6(w-2)}\geq\frac{q}{q-1}$. Since $w(w-1)=(w-2)(w+1)+2$, we
        have $6h(w)=w+1+\frac{2}{w-2}$, and hence
        $6\bigl(h(w+1)-h(w)\bigr)=1-\frac{2}{(w-1)(w-2)}>0$ for $w\geq4$; that
        is, $h$ is increasing for $w\geq4$. As $w\geq q+2\geq7$ in this subcase,
        it is enough to check the inequality at the smallest admissible value
        $w=q+2$, where it becomes $\frac{(q+2)(q+1)}{6q}\geq\frac{q}{q-1}$, that
        is, $(q+2)(q+1)(q-1)\geq6q^2$, or $q^3-4q^2-q-2\geq0$; this holds for
        $q\geq5$, since the left-hand side equals $18$ for $q=5$ and is
        increasing for $q\geq3$.
    \end{itemize}

    \emph{Case $s\geq m/2$.} By Lemma~\ref{lem:supersat}(i) and~\eqref{eq:d0},
    applied to $H$, whose number of edges is $e\coloneqq \binom m2-f=\lfloor
    m^2/4\rfloor+s$,
    \[
    d_0-s\ \geq\ g(s)\coloneqq \frac{e(4e-m^2)}{3m}-s .
    \]
    Since $\frac{de}{ds}=1$, we have $\frac{d}{ds}g(s)=\frac{8e-m^2}{3m}-1$, and
    as $e\geq\frac{m^2-1}{4}+\frac m2$ for $s\geq m/2$, this is at least
    $\frac{m^2+4m-2}{3m}-1=\frac{m^2+m-2}{3m}>0$. Hence $g$ is increasing, and
    since $s$ is an integer with $s\geq m/2$, it suffices to evaluate it at
    $s=\lceil m/2\rceil$. If $m$ is even, then $e=\frac{m(m+2)}{4}$, $4e-m^2=2m$
    and $g(s)=\frac{m(m+2)}{6}-\frac m2=\frac{m(m-1)}{6}$; if $m$ is odd, then
    $e=\frac{(m+1)^2}{4}$, $4e-m^2=2m+1$ and
    $g(s)=\frac{(m+1)^2(2m+1)}{12m}-\frac{m+1}{2}=\frac{(m+1)(m-1)(2m-1)}{12m}$.
    In both cases $g(s)\geq m-2\geq w-2$ for $m\geq4$ (for $m$ even this is
    $(m-3)(m-4)\geq0$, and for $m$ odd it is $2m^3-13m^2+22m+1\geq0$, which
    holds for $m\geq3$), so $d_0\geq s+w-2$, which gives~\eqref{eq:goal}. 
\end{proof}

\begin{figure}[htbp]
    \centering
    \includegraphics[width=0.9\textwidth]{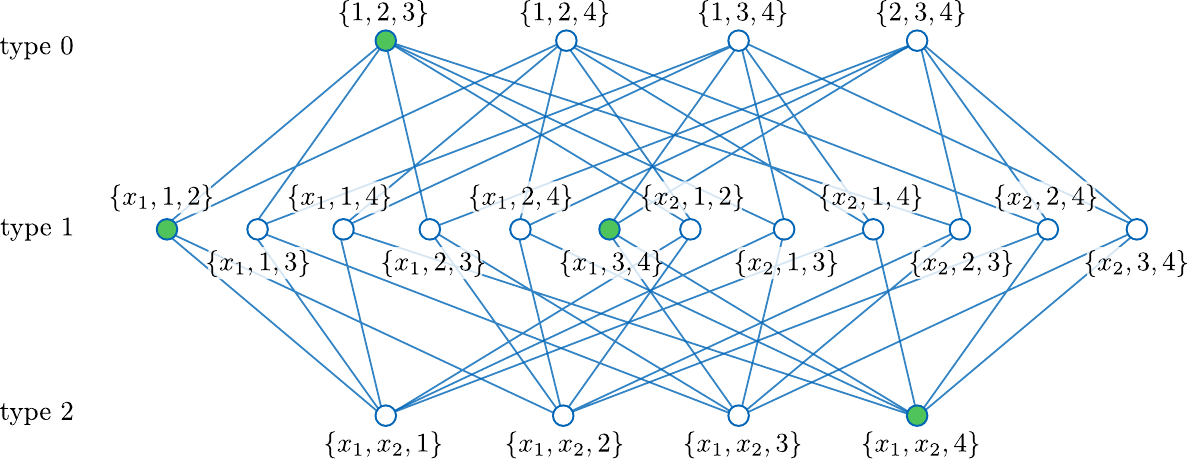}
    \caption{The graph $F_3(K_{2,4})$, with $V_1=\{x_1,x_2\}$ and
    $V_2=\{1,2,3,4\}$, drawn by types.  The filled vertices form the dominating
    set of the proof of Theorem~\ref{thm:K2n3} for $B_1=\{1,2\}$ and
    $B_2=\{3,4\}$: the configurations $\{x_1,1,2\}$ and $\{x_1,3,4\}$ given by
    the family~\eqref{eq:mantel}, the configuration $\{x_1,x_2,4\}$, and
    $T_0=\{1,2,3\}$. It has $T(4,3,2)+n-4=4$ elements, which is
    $\gamma(F_3(K_{2,4}))$.}
    \label{fig:K24}
\end{figure}

We close this section by studying the domination number of the 3-token graphs of
balanced complete bipartite graphs.  We first consider the small cases
$n\in\{6,8,10\}$.

\begin{remark}
Using integer linear programming, we obtain
\[
    \gamma(F_3(K_{3,3}))=2,\ \gamma(F_3(K_{4,4}))=6,\ \gamma(F_3(K_{5,5}))=12.
\]

We now give explicit minimum dominating sets attaining these values. Let
$V_1=\{x_1,\ldots,x_m\}$ and $V_2=\{y_1,\ldots,y_m\}$ be the partite sets of
$K_{m,m}$.

For $m=3$, consider
\[
D_6=\{\{x_1,x_2,x_3\},\{y_1,y_2,y_3\}\}.
\]

For $m=4$, let
\[
D_8=
\big\{\{ x_1,x_3,x_4\},\ \{y_2,y_3,y_4\},\ \{x_1,x_2,y_1\},\ \{x_3,x_4,y_1\},\ \{x_2,y_1,y_4\},\ \{x_2,y_2,y_3\}\big\}.
\]

Finally, for $m=5$, let
\[
\begin{aligned}
D_{10}=\big\{&
\{x_1,x_2,x_3\},\ \{x_1,x_2,x_4\},\ \{x_1,x_3,x_4\},\ \{y_1,y_2,y_4\},\ \{y_1,y_2,y_5\},\ \{y_2,y_4,y_5\},\\
&\{x_1,x_5,y_3\},\ \{x_2,x_4,y_3\},\ \{x_3,x_5,y_3\},\ \{x_5,y_1,y_5\},\ \{x_5,y_2,y_3\},\ \{x_5,y_3,y_4\}
\big\}.
\end{aligned}
\]
These sets are minimum dominating sets of
$F_3(K_{3,3})$, $F_3(K_{4,4})$, and
$F_3(K_{5,5})$, respectively.
\end{remark}

We now turn to the general case of balanced complete bipartite graphs with three tokens for even $n\geq 12$, using again a Turán-type argument.

\thmKmm*

\begin{proof}
    Let $m=n/2\geq6$, and let $V_1$ and $V_2$ be the partite sets of $K_{m,m}$.
    A configuration is said to be of type $(i,j)$ if it contains $i$ vertices of
    $V_1$ and $j$ vertices of $V_2$, where $i+j=3$.

    We first establish the lower bound. Let $D$ be a minimum dominating set of
    $F_3(K_{m,m})$.  By Proposition~\ref{prop_lwbnd_token} and Mantel's theorem,
    we have
    \[
    |D|\geq 2T(m,3,2)=2\binom{m}{2}
    -2\left\lfloor\frac{m^2}{4}\right\rfloor
    \geq \frac{m^2}{2}-m.
    \]
    Thus, if $m\leq 84$, then
    \[
    \frac{7}{48}n^2-4n=\frac{7}{12}m^2-8m
    \leq \frac{m^2}{2}-m,
    \]
    and the result follows.

    Now assume that $m\geq 85$, and note that if $|D|>\frac{7}{12}m^2$, the
    result follows. Hence we may also assume that
    \begin{equation}\label{eq_D_ass}
    |D|\leq\frac{7}{12}m^2. 
    \end{equation}

    For $i+j=3$, let $a_{i,j}$ denote the number of configurations of type
    $(i,j)$ belonging to $D$. We begin by establishing lower bounds on $a_{2,1}$
    and $a_{1,2}$.

    Let $\mathcal F_{2,1}$ be the family of pairs $P\subseteq V_1$ such that
    $P\cup\{b\}\in D$ for some $b\in V_2$.  Let $H_1$ be the graph on $V_1$
    whose edges are precisely the pairs not belonging to $\mathcal F_{2,1}$.
    Every triangle of $H_1$ corresponds to a configuration of type $(3,0)$ that
    cannot be dominated by a configuration of type $(2,1)$ belonging to $D$.
    Consequently, this configuration must itself belong to $D$, and therefore
    \[
    a_{3,0}\geq \operatorname{tr}(H_1). 
    \]
    Let $s\coloneqq T(m,3,2)-|\mathcal F_{2,1}|$. If $s\geq 1$, then
    \[
    |E(H_1)|
    =\left\lfloor\frac{m^2}{4}\right\rfloor+s.
    \]
    By Lemma \ref{lem:supersat} (i), we have
    \[
    \operatorname{tr}(H_1)
    \geq
    \frac{|E(H_1)|(4|E(H_1)|-m^2)}{3m}.
    \]
    For $m\geq85$, this implies $a_{3,0}\geq \frac{ms}{5}$.
    Thus, whether $s$ is positive or not, we have $s\leq \frac{5\,a_{3,0}}{m}$.
    Since $a_{2,1}\geq|\mathcal F_{2,1}|$, we obtain
        \[
    	a_{2,1}
    	\geq T(m,3,2)-s
    	\geq \frac{m^2}{4}-\frac m2
    	-\frac{5\,a_{3,0}}m
    	\geq \frac{m^2}{4}-\frac m2
    	-\left(\frac{5}{m}\right)\frac{7m^2}{12}
    	\geq \frac{m^2}{4}-4m.
    \]
    Analogously, by symmetry, we have
    \[
    a_{1,2}\geq\frac{m^2}{4}-4m.
    \]

    Now we focus on the domination of configurations of type $(2,1)$. For each
    $b\in V_2$, define
    \[
    U_b\mathrel{\coloneqq}\left\{
    u\in V_1:
    \{u,b,b'\}\in D
    \text{ for some }b'\in V_2\setminus\{b\}
    \right\}.
    \]
    Let $x_b\coloneqq m-|U_b|$, and let $e_b$ be the number of configurations of type
    $(2,1)$ in $D$ that contain $b$. Consider a pair $\{u,u'\}\subseteq
    V_1\setminus U_b$.  The configuration $\{u,u',b\}$ cannot be dominated by
    any configuration of type $(1,2)$ belonging to $D$. Hence either
    $\{u,u',b\}\in D$, or $\{u,u'\}$ is contained in a configuration of type
    $(3,0)$ belonging to $D$. Since each configuration of type $(3,0)$ contains
    exactly three pairs, it follows that
    \[
    \binom{x_b}{2}\leq 3\,a_{3,0}+e_b.
    \]
    By summing over $b\in V_2$ we obtain
    \begin{equation} \label{eq:double_domination}
    \sum_{b\in V_2}\binom{x_b}{2}
    \leq 3m\, a_{3,0}+a_{2,1}.
    \end{equation}

    On the other hand, each configuration $\{u,b,b'\}\in D$ of type $(1,2)$
    contributes $u$ to at most two sets $U_b$ and $U_{b'}$. Therefore,
    \[
    \sum_{b\in V_2}|U_b|\leq 2\,a_{1,2}.
    \]
    Consequently,
    \begin{equation}\label{eq:average_x}
    \frac1m\sum_{b\in V_2}x_b
    \geq m-\frac{2\,a_{1,2}}m.
    \end{equation}
    Now set 
    \[\overline{x}\mathrel{\coloneqq}m-\frac{2\,a_{1,2}}m, \] and write
    \[
    a_{2,1}=\frac{m^2}{4}+y_1,
    \qquad
    a_{1,2}=\frac{m^2}{4}+y_2.
    \]
    Hence $y_1,y_2\geq -4m$. If $y_2\geq\frac{m^2}{4}-\frac m4$,
    we obtain 
    \[
    |D|
    \geq a_{2,1}+a_{1,2}
    \geq \frac{3m^2}{4}-\frac{17m}{4}
    >\frac{7}{12}m^2,
    \]
    where the last inequality holds for $m\geq85$. This contradicts
    \eqref{eq_D_ass}. Thus we must have
    $y_2<\frac{m^2}{4}-\frac m4$,
    and hence
    \[
    \overline{x}
    =m-\frac{2\,a_{1,2}}m
    =\frac m2-\frac{2y_2}m
    >\frac12.
    \]
   Moreover, since $y_2\geq-4m$, we have
    \[
    \overline{x}\leq\frac m2+8.
    \]
    Now, since the function defined by
    \[
    f(x)=\binom{x}{2}=\frac{x(x-1)}2
    \]
    is convex and increasing for $x\geq1/2$, by Jensen's inequality and
    \eqref{eq:average_x} we have
       \[
    	\sum_{b\in V_2}\binom{x_b}{2}
    	\geq m\binom{\overline{x}}2
    	=\frac m2\overline{x}^{\,2}
    	-\frac m2\overline{x}
    	\geq
    	\frac m2\left(\frac m2-\frac{2y_2}m\right)^2
    	-\frac m2\left(\frac m2+8\right)
    	\geq
    	\frac{m^3}{8}-my_2-\frac{m^2}{4}-4m.
    \]
    Combining this inequality with \eqref{eq:double_domination}, we obtain
    \[
    3m\, a_{3,0}+a_{2,1}
    \geq
    \frac{m^3}{8}-my_2-\frac{m^2}{4}-4m.
    \]
    Using $a_{2,1}\leq |D|\leq7m^2/12$, we conclude that
    \[
    a_{3,0}
    \geq
    \frac1{3m}
    \left(
    \frac{m^3}{8}-my_2
    -\frac{m^2}{4}-4m-\frac{7m^2}{12}
    \right)
    \geq \frac{m^2}{24}-\frac{y_2}{3}-m.
    \]
    Analogously, by symmetry, we have that
    \[
    a_{0,3}\geq
    \frac{m^2}{24}-\frac{y_1}{3}-m.
    \]
    Finally, we have that
        \[
        \begin{aligned}
    	|D|&
    	=a_{3,0}+a_{2,1}+a_{1,2}+a_{0,3}
    	\geq
    	\frac{m^2}{12}
    	-\frac{y_1+y_2}{3}
    	-2m+\frac{m^2}{2}+y_1+y_2 \\
    	&=\frac{7m^2}{12}
    	+\frac23(y_1+y_2)-2m
    	\geq\frac{7m^2}{12}-\frac{22m}{3}
    	\geq\frac{7m^2}{12}-8m
    	 = \frac{7}{48}n^2-4n. 
    	 \end{aligned}
    \]
    This proves the lower bound.

    \medskip

    We now prove the upper bound. Let $V_1=A_1\sqcup A_2$ and $V_2=B_1\sqcup
    B_2$ be balanced bipartitions, with $|A_1|=|B_1|=\lceil m/2\rceil$ and
    $|A_2|=|B_2|=\lfloor
    m/2\rfloor$. We will use the following well-known fact (see, e.g.,
    \cite{Diestel}): for every $q\geq2$, the chromatic index of $K_q$ is at most
    $q$, that is, the edges of $K_q$ can be properly colored with $q$ colors. In
    other words, the pairs of a $q$-set can be labeled with $q$ labels so that
    the $q-1$ pairs containing any given element receive distinct labels.

    For $i\in\{1,2\}$, label the pairs of $\binom{A_i}{2}$ with labels in $B_1$
    in such a way that the pairs containing any given element of $A_i$ receive
    distinct labels; this is possible by the fact above, since $|B_1|\geq|A_i|$.
    Denote by $\lambda(P)\in B_1$ the label of a pair $P$, and let
    \[
        D_1  \coloneqq 
        \Bigl\{P\cup\{\lambda(P)\}:P\in\tbinom{A_1}{2}\cup\tbinom{A_2}{2}\Bigr\},
    \]
    so that 
    \[
    |D_1|=\binom{\lceil m/2\rceil}{2}+\binom{\lfloor
    m/2\rfloor}{2}=T(m,3,2).
    \]
    For every $u\in V_1$, the pairs containing $u$
    receive $|A_i|-1\geq|B_1|-2$ distinct labels, where $A_i$ is the class of
    $u$; for each of the at most two labels $b\in B_1$ not received by any pair
    containing $u$, add to a set $D_1'$ one configuration $\{u,u',b\}$, where
    $u'\in V_1\setminus\{u\}$ is arbitrary. Then $|D_1'|\leq 2m$, and the family
    $D_1\cup D_1'$ has the following property:
    \begin{equation}\label{eq:label}
        \text{for every $u\in V_1$ and every $b\in B_1$
        there is $u'\in V_1$ with $\{u,u',b\}\in D_1\cup D_1'$.}
    \end{equation}
    Symmetrically, label the pairs of $\binom{B_1}{2}\cup\binom{B_2}{2}$ with
    labels in $A_1$, and let $D_2$ and $D_2'$ be the corresponding families,
    with $|D_2|=T(m,3,2)$, $|D_2'|\leq2m$, and such that
    \begin{equation}\label{eq:label2}
        \text{for every $b\in V_2$ and every $a\in A_1$
        there is $b'\in V_2$ with $\{a,b,b'\}\in D_2\cup D_2'$.}
    \end{equation}
    Finally, let $D_3\subseteq\binom{B_2}{3}$ be a covering of the pairs of
    $B_2$ by triples of size $C(\lfloor m/2\rfloor,3,2)$, and let
    $D_4\subseteq\binom{A_2}{3}$ be a covering of the pairs of $A_2$ of the same
    size. By~\eqref{eq:FH}, we have 
    \[
    |D_3|=|D_4|=\frac{m^2}{24}+O(m).
    \]
    Let
    \[
        D \coloneqq D_1\cup D_1'\cup D_2\cup D_2'\cup D_3\cup D_4.
    \]
    We claim that $D$ is a dominating set of $F_3(K_{m,m})$. To show this, let
    $A$ be a configuration. Now we consider the following cases.

    \begin{itemize}
        \item If $A$ is of type $(3,0)$, then two of its elements lie in the
        same class $A_i$, so $A$ contains a pair $P\in\binom{A_i}{2}$, and
        $A\sim P\cup\{\lambda(P)\}\in D_1$, since their symmetric difference
        consists of the element of $A\setminus P$ and $\lambda(P)\in V_2$.
        Symmetrically, if $A$ is of type $(0,3)$ then $A$ is dominated by $D_2$.
        
        \item If $A=\{u,b,b'\}$ is of type $(1,2)$ and $\{b,b'\}\cap
        B_1\neq\varnothing$, say $b\in B_1$, then by~\eqref{eq:label} there is
        $\{u,u',b\}\in D_1\cup D_1'$, and $A\sim\{u,u',b\}$, since their
        symmetric difference is $\{b',u'\}$. If $\{b,b'\}\subseteq B_2$, let
        $T\in D_3$ be a triple containing $\{b,b'\}$; then $A\sim T$, since
        their symmetric difference is $\{u,b''\}$ with $b''\in V_2$.

        \item Symmetrically, if $A=\{u,u',b\}$ is of type $(2,1)$, then $A$ is
        dominated by $D_2\cup D_2'$ if $\{u,u'\}\cap A_1\neq\varnothing$,
        by~\eqref{eq:label2}, and by $D_4$ if $\{u,u'\}\subseteq A_2$.
    \end{itemize}
    Hence $D$ is a dominating set, and
    \[
        |D| \leq 2\,T(m,3,2)+2\,C\left(\left\lfloor\tfrac m2\right\rfloor,3,2\right)+4m
        = \frac{m^2}{2}+\frac{m^2}{12}+O(m)=\frac{7}{12}m^2+O(m).
    \]
    Since $m=n/2$, it follows that $|D|\leq\frac{7}{48}n^2+O(n)$, which
    completes the proof.
\end{proof}

\section*{Acknowledgments}

This work began during the 6th Workshop on Optimization, Mathematics and
Algorithms
(ROMA~2025)\footnote{\url{https://academicos.azc.uam.mx/heredia/ROMA/2025/}}. We
thank the participants for fruitful scientific discussions.  G.
Ben\'itez-Bobadilla is supported by postdoctoral grant “Estancias Posdoctorales
por M\'exico” by SECIHTI (CVU: 622815).  F.~E.~Contreras-Mendoza gratefully
acknowledges support from the UNAM-PAPIIT IN106425 grant.

\section*{Declaration of AI use}

During the preparation of this work, the authors used Claude (Anthropic) as an
assistant for literature search, for exploring and checking arguments, for
computational verification of small cases, and for drafting and editing parts of
the text. All results, proofs and computations were reviewed and verified by the
authors, who take full responsibility for the content of this paper. 

\bibliographystyle{abbrv}
\bibliography{Domination}

\end{document}